\documentclass[3p,times]{elsarticle}

\usepackage{ecrc}
\usepackage{algorithm} 
\usepackage{multirow} 
\usepackage{amsmath}
\usepackage{xcolor}
\usepackage{amsmath, amssymb, epsfig, graphicx, bm}
\usepackage{cases}
\usepackage{float, subfigure}
\usepackage{color}
\usepackage{epstopdf}
\usepackage{mathrsfs}
\usepackage{amsthm}
\usepackage{amsfonts}
\usepackage{url}
\usepackage[all,cmtip]{xy}
\usepackage{color,hyperref}
\usepackage[noend]{algpseudocode}
\usepackage{algorithmicx,algorithm}
\usepackage{setspace}
\usepackage[misc]{ifsym}

\newcommand{\black}[1]{{\color{black}#1}}
\newcommand{\blue}[1]{{\color{blue}#1}}

\definecolor{darkblue}{rgb}{0,0,1}
\hypersetup{colorlinks,breaklinks,
	linkcolor=darkblue,urlcolor=darkblue,anchorcolor=darkblue,citecolor=darkblue}

\usepackage[normalem]{ulem} 

\newtheorem{theorem}{Theorem}

\newtheorem{assumption}{Assumption}

\newcommand{\R}{{\mathbb{R}}}

\newcommand{\E}{\mathbb{E}}

\volume{00}

\firstpage{1}

\journalname{Signal Processing}

\runauth{}

\jid{procs}

\jnltitlelogo{SIGPRO}

\usepackage{amssymb}

\usepackage[figuresright]{rotating}

\begin{document}

	\begin{frontmatter}
		
		
		
		\dochead{}
		
        \title{Byzantine-Robust Decentralized Stochastic Optimization \\ over Static and Time-Varying Networks}


        \author{Jie Peng \quad Weiyu Li \quad  Qing Ling\textsuperscript{*} }

        \address{  \textsuperscript{*}Corresponding Author \\ \ E-mail address: lingqing556@mail.sysu.edu.cn \\ School of Data and Computer Science and Guangdong Province Key Laboratory of Computational Science, \\ Sun Yat-Sen University, Guangzhou, Guangdong, 510006, China}

        \address{ }

		\begin{abstract}
			In this paper, we consider the Byzantine-robust stochastic optimization problem defined over decentralized static and time-varying networks, where the agents collaboratively minimize the summation of expectations of stochastic local cost functions, but some of the agents are unreliable due to data corruptions, equipment failures or cyber-attacks. The unreliable agents, which are called as Byzantine agents thereafter, can send faulty values to their neighbors and bias the optimization process. Our key idea to handle the Byzantine attacks is to formulate a total variation (TV) norm-penalized approximation of the Byzantine-free problem, where the penalty term forces the local models of regular agents to be close, but also allows the existence of outliers from the Byzantine agents. A stochastic subgradient method is applied to solve the penalized problem. We prove that the proposed method reaches a neighborhood of the Byzantine-free optimal solution, and the size of neighborhood is determined by the number of Byzantine agents and the network topology. Numerical experiments corroborate the theoretical analysis, as well as demonstrate the robustness of the proposed method to Byzantine attacks and its superior performance comparing to existing methods.
		\end{abstract}
		
		\begin{keyword}
			Decentralized stochastic optimization, Byzantine attacks, robustness, static networks, time-varying networks
		\end{keyword}
		
	\end{frontmatter}
	
	
	
	\section{Introduction}
	\label{sec:1}
	In recent years, decentralized stochastic optimization has become a popular research topic in the signal processing and machine learning communities. With the rapidly increasing number of distributed devices and volume of generated data, traditional signal processing and machine learning approaches, which rely on a central controller to collect the data samples or coordinate the optimization process, suffer from privacy and scalability issues \cite{KoloskovaSJ19}. In decentralized stochastic optimization, every device (called as agent thereafter) learns its own model using its local data samples, and periodically exchanges its model with neighboring agents so as to achieve consensus. This scheme is favorable in privacy preservation since the data samples are kept local, and does not rely on any central controller that could be a system bottleneck. Existing decentralized stochastic optimization methods include decentralized parallel stochastic gradient descent (DPSGD) \cite{dpsgd}, stochastic subgradient projection \cite{sspm}, dual averaging \cite{dam}, mirror descent \cite{mdm}, etc. Asynchronous algorithms are developed in \cite{adpsgd,aapsgd} to reduce the idle time, and variance reduction techniques are proposed in \cite{d2,vrsl,dsa,chi-vr,gtvr-survey} to improve the convergence rate. Decentralized stochastic optimization methods are shown to be superior to their centralized counterparts on training large-scale neural networks when the communication links are subject to high latency and limited bandwidth \cite{hop}.
	
	However, the lack of centralized coordination in decentralized stochastic optimization also raises concerns on robustness. Some of the agents might be malfunctioning or even malicious. Due to data corruptions, equipment failures or cyber-attacks, they can send faulty values to their neighbors and bias the optimization process. \black{We consider a general Byzantine attack model \cite{LamportSP82}, in which the Byzantine agents are omniscient and can arbitrarily modify the values sent to other agents. Such an attack model imposes no restrictions on the adversaries and is worst-case.} The purpose of this paper is to develop a Byzantine-robust decentralized stochastic optimization method.
	
	Most of the existing decentralized stochastic optimization methods are vulnerable to Byzantine attacks. Take DPSGD as an example. At every iteration, every agent mixes the models received from its neighbors, followed by a stochastic gradient step on the cost function constructed from one local data sample (or a batch of them), to update its local model \cite{dpsgd}. When the Byzantine agents send well-designed faulty values instead of the true models, they are able to lead the regular agents to end up with incorrect results.
	
	Byzantine-robust \textit{decentralized deterministic} optimization methods have been developed in \cite{ByRDiE,BRIDGE}, where at every iteration every regular agent uses all of its local data samples, instead of one or a batch. The work of \cite{ByRDiE} proposes ByRDiE, in which every regular agent utilizes coordinate-wise trimmed mean to screen outliers in the received models, and then applies coordinate gradient descent to update its local model. The one-coordinate-at-a-time update of ByRDiE is inefficient for high-dimensional problems \cite{BRIDGE}. To address this issue, the work of \cite{BRIDGE} proposes BRIDGE, which allows every regular agent to update all the coordinates of its local model at every iteration. Although these two algorithms are originally developed for \textit{decentralized deterministic} optimization, they can also be adapted to the \textit{decentralized stochastic} setting according to our numerical experiments. However, to the best of our knowledge, most of the existing works do not explicitly consider the Byzantine-robust \textit{decentralized stochastic} optimization problem; see for reference the recent survey paper \cite{adversary-resilient}.
	
	
	There are some works that consider the Byzantine-robust \textit{centralized stochastic} optimization problem, where a central controller aggregates the information from the agents and coordinates the optimization process. The main idea of these works is to modify the stochastic gradient method with robust aggregation rules. To be specific, at every iteration, the central controller sends the current model to all the agents, the regular agents send back their local stochastic gradients, while the Byzantine agents may send back faulty values. When the local data samples are independently and identically distributed (i.i.d.), the local stochastic gradients are also i.i.d. and the central controller can obtain a reliable approximation to the average of the local stochastic gradients through aggregating all the received values with trimmed mean, geometric median, or other robust aggregation rules \cite{Chen2019DistributedSM,Xie2018GeneralizedBS}. However, this idea is not directly applicable to \textit{decentralized stochastic} optimization. Since there is no central controller to maintain a common model, the regular agents have to evaluate their local stochastic gradients at different points. Therefore, even though the local data samples are i.i.d. the local stochastic gradients are not necessarily so, and thus the robust aggregation rules have no theoretical guarantee in this case.
	
	In addition to Byzantine-robust decentralized stochastic optimization over static networks, we are also interested in the case that the underlying network topologies are \textit{time-varying}, which finds applications in mobile computation systems and robotic swarms. Various decentralized stochastic optimization algorithms over time-varying networks have been developed, such as those in \cite{yuan2014randomized, nedic2017achieving}. Some of the recent works further consider that the time-varying networks are directed such that the designed mixing matrices are column stochastic \cite{nedic2016stochastic,nedic2014distributed,scutari2019distributed} or row stochastic \cite{hong2017stochastic}. However, none of these algorithms take into account the existence of Byzantine agents. \black{On the other hand, there exist other works considering adversarial attacks over time-varying graphs. For example, the works of \cite{chen2019resilient, chen2020resilient} focus on robust distributed parameter estimation, while their attack model is different from ours. To be specific, the works of \cite{chen2019resilient, chen2020resilient} consider adversarial attacks on the measured data, while our Byzantine attack model considers attacks during the optimization process.}

	This paper develops a Byzantine-robust \textit{decentralized stochastic} optimization method over static and time-varying networks, where the networks are fully decentralized and contain unknown numbers of Byzantine agents, the local data samples at the regular agents are not necessarily i.i.d. and only one data sample (or a batch of them) is available for every regular agent at every iteration. The key idea is to formulate a total variation (TV) norm-penalized approximation of the Byzantine-free problem, where the penalty term forces the local models of regular agents to be close, but also allows the existence of outliers from the Byzantine agents. A stochastic subgradient method is applied to solve the penalized problem (Section \ref{sec:2}). Although the TV norm-penalized approximation has been investigated in Byzantine-robust \textit{decentralized deterministic} \cite{sta-opt}, \textit{decentralized dynamic} \cite{dyn-opt} and \textit{centralized stochastic} \cite{RSA} optimization problems, its application in Byzantine-robust \textit{decentralized stochastic} optimization is novel. We prove that the proposed method reaches a neighborhood of the Byzantine-free optimal solution under mild assumptions, and the size of neighborhood is determined by the number of Byzantine agents and the network topology (Section \ref{sec:3}). Numerical experiments corroborate the theoretical analysis and demonstrate the robustness of the proposed method to Byzantine attacks (Section \ref{sec:4}).
	
	
	\section{Problem Statement and Algorithm Development}
	\label{sec:2}
	
	In this section, we describe the formulation of Byzantine-robust decentralized stochastic optimization problem, and develop the algorithm for both static and time-varying networks.
	
	\subsection{Static Network Case}
	\label{static-network}
	Consider a static and undirected network $\mathcal{G} = (\mathcal{V}, \mathcal{E})$ with a set of $n$ agents $\mathcal{V}= \{1,\cdots,n\}$ and a set of undirected edges $\mathcal{E}$. If $e = (i,j) \in \mathcal{E}$, then agents $i$ and $j$ are neighbors and can communicate with each other at a low cost. Since the network is undirected, for notational convenience, we let every $e = (i, j)$ satisfy $i < j$. However, not all the agents are regular. An unknown number of Byzantine agents are supposed to be omniscient and can send faulty values to their neighbors during the optimization process. Denote $\mathcal{R}$ and $\mathcal{B}$ as the sets of regular agents and Byzantine agents, respectively. We have $\mathcal{V} = \mathcal{R} \cup \mathcal{B}$. For agent $i$, denote the set of its regular neighbors as $\mathcal{R}_i$ and the set of its Byzantine neighbors as $\mathcal{B}_i$. Thus, $\mathcal{N}_i := \mathcal{R}_i \cup \mathcal{B}_i$ is the set of all neighbors of agent $i$. Denote $\mathcal{E}_{R} \subseteq \mathcal{E}$ as the set of reliable edges not attached to any Byzantine agent. The decentralized stochastic optimization problem defined over the network is
	\begin{align}\label{eq4}
		\tilde{x}^* = \arg\min\limits_{\tilde{x} \in \R^p} \sum_{i \in \mathcal{R}} \bigg(\E[F(\tilde{x},\xi_i)] + f_0(\tilde{x})\bigg),
	\end{align}
	where $\tilde{x} \in \R^p$ is an optimization variable (also called as model), $F(\tilde{x},\xi_i)$ is a smooth cost function determined by a random variable $\xi_i$ following a distribution $\mathcal{D}_i$ and represents an empirical loss related to a randomly chosen data sample in regular agent $i$, and $f_0(\tilde{x})$ is a smooth regularization term. \black{Instead of the i.i.d. assumption as in \cite{ByRDiE, BRIDGE}, here the random variables $\{\xi_i, i \in \mathcal{R}\}$ are only assumed to be mutually independent. We will also assume $\E[F(\tilde{x}, \xi_i)] + f_0(\tilde{x})$ to be strongly convex in Assumption \ref{assumption1} in Section \ref{sec:3}.} Our goal is to find the optimal solution $\tilde{x}^*$ through collaboration of the regular agents. The main challenges are three-fold: (i) the network lacks a central coordinator and is fully decentralized, (ii) only one randomly chosen data sample (or a batch of them) can be used by every regular agent at every iteration, and, (iii) more importantly, the Byzantine agents can send faulty values to their neighbors so as to bias the optimization process, but their identities are unknown.
	
	To develop a reasonable algorithm, it is necessary to assume that the network of regular agents is  connected \cite{dyn-opt}. Otherwise, if a regular agent is surrounded by Byzantine neighbors, it is unable to communicate and collaborate with any regular agents. Therefore, the best model it can learn is solely based on its local data samples, and may be far away from the true model in the non-i.i.d. setting.
	\begin{assumption} \label{assumption-static}
		\textbf{(Network Connectivity over Static Graph)} The network consisting of all regular agents $i \in \mathcal{R}$, denoted as $(\mathcal{R},\mathcal{E}_{R})$, is connected.
	\end{assumption}

	We begin from assuming that the Byzantine agents are absent. Rewrite \eqref{eq4} to a consensus-constrained form, which is common in decentralized optimization. Denote $x_i \in \R^p$ as the local copy of the model $\tilde{x}$ at regular agent $i$ and stack all local copies in a longer vector $x:=[x_i]\in \R^{|\mathcal{R}|p}$. When the regular agents are connected as stated by Assumption \ref{assumption-static}, \eqref{eq4} is equivalent to
	\begin{align}\label{eq5}
		\min\limits_{x:=[x_i]} & ~ \sum_{i \in \mathcal{R}} \bigg(\E[F(x_i,\xi_i)] + f_0(x_i)\bigg), \\
		s.t.                   & ~ x_i = x_j, ~ \forall i \in \mathcal{R}, ~ \forall j \in \mathcal{R}_i, \nonumber
	\end{align}
	in the sense that $[\tilde{x}^*] \in \R^{|\mathcal{R}|p}$ that stacks $|\mathcal{R}|$ vectors of $\tilde{x}^*$, the optimal solution to \eqref{eq4}, is the optimal solution to \eqref{eq5}.
	
	Then, motivated by \cite{sta-opt,dyn-opt,RSA}, we propose to solve a TV norm-penalized approximation of \eqref{eq5}, as
	\begin{align}\label{eq6}
		x^* = \arg\min\limits_{x:=[x_i]} ~ \sum_{i \in \mathcal{R}} \bigg( \E[F(x_i, \xi_i)] + \frac{\lambda}{2} \sum_{j \in \mathcal{R}_i} \|x_i  -  x_j \|_1 + f_0(x_i)\bigg),
	\end{align}
	where $\lambda \geq 0$ is a penalty parameter. For every pair of regular neighbors $(i,j)$, $x_i$ and $x_j$ are forced to be close through introducing the TV norm penalty \black{$\frac{\lambda}{2} \sum_{i \in \mathcal{R}} \sum_{j \in \mathcal{R}_i} \|x_i  -  x_j \|_1$, in which every pair of regular neighbors $(i,j)$ appears once}. The larger $\lambda$ is, the closer $x_i$ and $x_j$ are forced to be. On the other hand, the TV norm penalty also allows some pairs of $x_i$ and $x_j$ to be different, which is important when the Byzantine agents are present as we will discuss later.
	
	Since calculating the full subgradient of the cost function in \eqref{eq6} is time-consuming or even impossible, we solve \eqref{eq6} with the stochastic subgradient method. At time $k$, every regular agent $i$ updates its local model $x_i^{k+1}$ as
	\begin{align}\label{eq7}
		x^{k+1}_i  = x^k_i - \alpha^{k} \Bigg(\nabla  F(x^k_i,\xi_i^k) + \lambda \sum_{j \in \mathcal{R}_i} sign(x^k_i-x^k_j) + \nabla f_0(x_i^k) \Bigg),
	\end{align}
	where \black{$\xi_i^k \sim \mathcal{D}_i$} corresponds to the random data sample chosen \black{independently} by agent $i$ at time $k$, $sign(\cdot)$ is the element-wise sign function, and $\alpha^{k} > 0$ is a step size. Given $\beta \in \R$, $sign(\beta)$ equals to $1$ when $\beta > 0$, $-1$ when $\beta < 0$, and an arbitrary value within $[-1, 1]$ when $\beta = 0$. Observe that \eqref{eq7} is fully decentralized. To update $x^{k+1}_i$, a regular agent $i$ needs to evaluate its own local stochastic gradient $\nabla  F(x^k_i,\xi_i^k)$ and gradient $\nabla f_0(x_i^k)$, as well as combine the models $\{x^k_j, j \in \mathcal{R}_i\}$ received from its regular neighbors $\{j, j \in \mathcal{R}_i\}$.
	
	
	Now we consider how \eqref{eq7} performs when the Byzantine agents are present. \black{A Byzantine agent $j$ will not send its true model to its neighbors at time $k$. Instead, it sends an arbitrary vector $z_j^k \in \R^p$.} In this case, for a regular agent $i$, \eqref{eq7} becomes
	\begin{align}\label{eq:subgradient}
		x^{k+1}_i  = x^k_i - \alpha^{k} \Bigg(&\nabla  F(x^k_i,\xi_i^k)  +  \lambda   \sum_{j \in \mathcal{R}_i} sign(x^k_i-x^k_j) + \lambda\sum_{j \in \mathcal{B}_i} sign(x^k_i-z^k_j) + \nabla f_0(x_i^k) \Bigg).
	\end{align}

	The resulting Byzantine-robust decentralized stochastic optimization method is outlined in Algorithm 1. In \eqref{eq:subgradient}, observe that the elements of $sign(x^k_i-z^k_j)$ are in the range of $[-1,1]$, such that the influence of the faulty vector $z^k_j$ is limited, although $z^k_j$ can be arbitrary. We will theoretically justify the robustness of the proposed algorithm to Byzantine attacks in the subsequent section.
	
	Note that the TV norm penalty introduced here is based on the $\ell_1$ norm. Other norms such as $\ell_2$ and $\ell_\infty$ are also applicable, as recommended in \cite{RSA} for Byzantine-robust centralized stochastic optimization. We leave their development and analysis to our future work.
	
	\begin{algorithm}
		\caption{Byzantine-robust decentralized stochastic optimization over static graph} {\bf Input:} $x_i^0 \in \R^p$ for
		$i \in \mathcal{R}$, $\lambda >0$, and $\{\alpha^k, k = 0, 1, \cdots\}$.
		\begin{algorithmic}[1]
			\For{$k = 0, 1, \cdots$, every regular agent $i \in \mathcal{R}$}
			\State Broadcast its current model $x_i^k$ to all the neighbors.
			\State Receive $x^k_j$ from regular
			neighbors $j \in \mathcal{R}_i$ and $z^k_j$ from Byzantine neighbors $j \in
			\mathcal{B}_i$.
			\State Update local iterate $x^{k+1}_i$ according to
			\eqref{eq:subgradient}.
			\EndFor \State \textbf{end for}
		\end{algorithmic}
	\end{algorithm}
	
	
	
	\subsection{Time-Varying Network Case}
	\label{time-vary}

	We further consider the more challenging scenario that the network has time-varying communication edges. \black{As in the static case, we denote $\mathcal{R}$ and $\mathcal{B}$ as the sets of regular agents and Byzantine agents, respectively, such that $\mathcal{V} = \mathcal{R} \cup \mathcal{B}$. Given the fixed agents, there are a finite number of possible undirected graphs with different edges. We encode these candidate graphs by the edge sets $\mathcal{E}(\zeta)$ with $ \zeta =[\zeta_e]_{e=(i,j)\in\mathcal{V}\times\mathcal{V}}$, in which each entry $\zeta_e\in
		\{0,1\}$ of $\zeta$ means $e\in\mathcal{E}$ if $\zeta_e = 1$, and $e\notin\mathcal{E}$ otherwise. At time $k$, the current graph can be represented as $\mathcal{G}^k = \big(\mathcal{V}, \mathcal{E}(\zeta^k)\big)$, where $\zeta^k$ is a random vector, not necessarily independent with each other across time, but independent with $\xi_i^{k'}$ for all agent $i$ and time $k'$. Then we obtain a countable infinite random sequence $(\zeta^1,\zeta^2,\ldots,\zeta^k,\ldots)$ which completely characterizes the time-varying network during the optimization process. To simplify the notations, we use $\mathcal{E}^k = \mathcal{E}(\zeta^k)$ as an abbreviation, and the stochasticity of $\zeta^k$ will be implicitly included in the following discussions.} 
	For agent $i$ at time $k$, denote $\mathcal{R}_i^k$ as the set of its instantaneous regular neighbors, $\mathcal{B}_i^k$ the set of its instantaneous Byzantine neighbors, and $\mathcal{N}_i^k := \mathcal{R}_i^k \cup \mathcal{B}_i^k$ the set of all instantaneous neighbors. Denote $\mathcal{E}_{R}^{k} \subseteq \mathcal{E}^k$ as the set of reliable edges not attached to any Byzantine agent. 
	
	\black{Though our goal is still to solve the decentralized stochastic optimization problem defined in \eqref{eq4}, the time-varying network topology is an essential issue compared to the static case. To guarantee sufficiently frequent information exchange among all regular agents, we consider the \textit{average network} consisting of all regular agents and the edges that ergodically appear many times.} To be specific, 
	define $\bar{a}_e$ as the limiting 
	frequency of a reliable edge $e=(i, j)$ appearing, where $i, j \in \mathcal{R}$ and $i < j$, given by $$\bar{a}_e = \liminf_{K \rightarrow \infty} \frac{1}{K+1} \sum_{k=0}^K \zeta_e^k.$$
		Denote $\bar{\mathcal{E}}_R = \{e: \bar{a}_e > 0\}$, and then call $(\mathcal{R},\bar{\mathcal{E}}_{R})$ as the average network.
The following assumption guarantees that the average network is connected.
	
	\black{
		\begin{assumption} \label{assumption-timevary}
			\textbf{(Network Connectivity over Time-Varying Graph)} The empirical distribution of $ \{\zeta^{k}\}_{k=0}^K $ converges to the distribution of some random variable $\bar\zeta=[\bar{\zeta}_e]$ as $K\rightarrow\infty$, where the expectation of $\bar\zeta_e$  is exactly $\bar{a}_e$ at each regular edge $e$. The average network $(\mathcal{R},\bar{\mathcal{E}}_{R})$ is connected.
	\end{assumption}} 
	%
	
	Assumption \ref{assumption-timevary} is much weaker than Assumption \ref{assumption-static} in which there exists a path between any two regular agents at
all times. In contrast, Assumption \ref{assumption-timevary} allows the path to be temporarily disconnected, while ensuring the regular agents to exchange information in a sufficiently frequent manner.

	\hspace{0.5em}
	
	\noindent \textbf{Remark 1}. 
	Our time-varying network model is able to describe several common scenarios.
	
	\noindent (i) Randomly activated edges. At every time $k$, every reliable edge $e$ is connected with probability $p_e \in (0, 1]$. In this scenario, $\bar{a}_e = p_e$.
	
	\noindent (ii) Periodical network. From time $0$ to $T-1$, every reliable edge $e$ appears $t_e$ times. These $T$ network topologies reappear at the following times, with a period of $T$. In this scenario, $\bar{a}_e = \frac{t_e}{T}$.
	
	\noindent (iii) Quasi-periodical network. For any time span with length $T$, every reliable edge $e$ appears at least $t_e$ time. In this scenario, $\bar{a}_e \geq \frac{t_e}{T}$.
	
	\hspace{0.5em}
	
	Similar to the equivalent transformation from \eqref{eq4} to \eqref{eq5}, if the regular agents are connected \black{in the average network} as stated by Assumption \ref{assumption-timevary}, \eqref{eq4} is equivalent to
	\black{
		\begin{align}\label{eq5-time}
			\min\limits_{x:=[x_i]} & ~ \sum_{i \in \mathcal{R}} \bigg(\E_{\xi_i}[F(x_i,\xi_i)] + f_0(x_i)\bigg), \\
			s.t.                   & ~ x_i = x_j, ~ \forall (i,j) \in \bar{\mathcal{E}}_{R}. \nonumber
		\end{align}
	} 
	The TV norm-penalized approximation of \eqref{eq5-time} is given by
	\black{
		\begin{align}
			x^* &= \arg\min\limits_{x:=[x_i]} ~ \sum_{i \in \mathcal{R}} \bigg( \E_{\xi_i}[F(x_i, \xi_i)] + f_0(x_i)\bigg) + \lambda\sum_{e=(i,j)\in\bar{\mathcal{E}}_R}\bigg(\bar{a}_e \|x_i  -  x_j \|_1\bigg) \notag\\
			&= \arg \min_{x:=[x_i]} ~ \sum_{i \in \mathcal{R}} \Bigg(\E_{\xi_i}[F(x_i, \xi_i)]  + \frac{\lambda}{2} \ \E_{\bar\zeta}[\sum_{j \in \mathcal{R}_i(\bar\zeta)} \|x_i - x_j\|_1] + f_0(x_i) \Bigg),
			\label{tvn-eq2}
		\end{align}
		where $\mathcal{R}_i(\bar\zeta)$ is the set of regular neighbors of agent $i$ in terms of the average network $\big(\mathcal{R},\mathcal{E}(\bar\zeta)\big)$.}
		
		
		Therefore, akin to the static case, we solve \eqref{tvn-eq2} by the stochastic subgradient method as follows. At time $k+1$, every regular agent $i$ updates its local model as 
	%
	%
	%
	\begin{align}\label{tvn-eq4}
		x_i^{k+1} = x_i^k - \alpha^k \bigg(\nabla F(x_i^k, \xi_i^k) + \lambda\sum_{j \in \mathcal{R}_i^k} sign(x_i^k- x_j^k)  + \nabla f_0(x_i^k) \bigg).
	\end{align}
	
	\black{Now we turn to the case with Byzantine attacks.} 
	A Byzantine agent $j$ will send an arbitrary vector $z_j^k \in \R^p$ instead of the true model $x_j^k$ to its neighbors at time $k$. Therefore, in the presence of Byzantine agents, \eqref{tvn-eq4} becomes
	\begin{align}\label{tvn-eq6}
		x_i^{k+1} = x_i^k - \alpha^k \bigg(\nabla F(x_i^k, \xi_i^k) + \lambda\sum_{j \in \mathcal{R}_i^k} sign(x_i^k- x_j^k) + \lambda\sum_{j \in \mathcal{B}_i^k} sign(x_i^k- z_j^k) + \nabla f_0(x_i^k) \bigg).
	\end{align}

	The resulting Byzantine-robust decentralized stochastic optimization method over the time-varying network is outlined in Algorithm 2. With particular note, our proposed method has a consistent form for both static and time-varying networks. No specific design of any mixing matrix is needed to adapt to the time-variance.
	
	\begin{algorithm}
		\caption{Byzantine-robust decentralized stochastic optimization over time-varying graph} {\bf Input:} $x_i^0 \in \R^p$ for
		$i \in \mathcal{R}$, $\lambda >0$ and $\{\alpha^k, k = 0, 1, \cdots\}$.
		\begin{algorithmic}[1]
			\For{$k = 0, 1, \cdots$, every regular agent $i \in \mathcal{R}$}
			\State Broadcast its current model $x_i^k$ to all the neighbors.
			\State Receive $x_j^k$ from regular neighbors $j \in \mathcal{R}_i^k$ and $z^k_j$ from Byzantine neighbors $j \in \mathcal{B}_i^k$.
			\State Update local iterate $x^{k+1}_i$ according to 	\eqref{tvn-eq6}.
			\EndFor \State \textbf{end for}
		\end{algorithmic}
	\end{algorithm}
	
	\section{Performance Analysis}
	\label{sec:3}
	
	In this section, we theoretically analyze  the proposed Byzantine-robust decentralized stochastic optimization method in terms of convergence and robustness. We make the following assumptions, which are common in analyzing decentralized stochastic optimization methods.
	
	\begin{assumption} \label{assumption1}
		\textbf{(Strong Convexity)} \black{For any model $\tilde{x} \in \R^p$ and every regular agent $i \in \mathcal{R}$, $\E [F(\tilde{x}, \xi_i)] + f_0(\tilde{x})$ is strongly convex with constants $u_i$, where $\E [F(\tilde{x},\xi_i)]$ is the local cost function and $f_0(\tilde{x})$ is the regularization term.}
	\end{assumption}
	
	\begin{assumption} \label{assumption2}
		\textbf{(Lipschitz Continuous Gradients)} \black{For any model $\tilde{x} \in \R^p$ and every regular agent $i \in \mathcal{R}$, $\E [F(\tilde{x}, \xi_i)] + f_0(\tilde{x})$ is differentiable and has Lipschitz continuous gradients with constants $L_i$, where $\E [F(\tilde{x},\xi_i)]$ is the local cost function and $f_0(\tilde{x})$ is the regularization term.}
	\end{assumption}
	
	
	
	\begin{assumption} \label{assumption3}
		\textbf{(Bounded Variance)} Every regular worker $i \in \mathcal{R}$ samples data with random variables $\xi_i^k \sim \mathcal{D}_i$ at every time $k$ \black{independently}. \black{For any model $\tilde{x} \in \R^p$, }the variance of $\nabla F(\tilde{x},\xi_i^k)$ is upper bounded by $\delta_i^2$, i.e., $\black{\E_{\xi_i^k}}[\| \nabla F(\tilde{x},\xi_i^k) - \black{\E_{\xi_i^k}}[\nabla F(\tilde{x},\xi_i^k)] \| ^2] \leq \delta_i^2, \forall i \in \mathcal{R}$.
	\end{assumption}
	
	The analysis in this paper shares similarities with that in \cite{RSA}. However, \cite{RSA} considers Byzantine-robust centralized stochastic optimization, while this paper considers the decentralized case. Due to the underlying decentralized static and time-varying networks, our proofs are significantly different from those in \cite{RSA}. Our theoretical results also explicitly show the influence of the topologies on the performance. Due to the page limit, we delegate the detailed proofs to \cite{Peng2020arxiv}.
	
	\subsection{Performance Analysis over Static Network}
	
	The first theorem shows that the TV norm-penalized problem \eqref{eq6} is equivalent to the consensus-constrained one \eqref{eq5} (and hence \eqref{eq4} too), when the penalty parameter $\lambda$ is sufficiently large. This theorem is analogous to Theorem 1 in \cite{RSA}, but our proof is based on a system of linear equations involving the decentralized network structure and different from the proof of Theorem 1 in \cite{RSA}.
	We define $A \in \mathbb{R}^{|\mathcal{R}| \times |\mathcal{E}_{R}|}$ as the node-edge incidence matrix of $(\mathcal{R},\mathcal{E}_{R})$. To be specific, for an edge $e=(i,j) \in \mathcal{E}_{R}$ with $i<j$, the $(i,e)$-th entry of $A$ is $1$ while the $(j,e)$-th entry of $A$ is $-1$.
	\begin{theorem} \label{theorem1}
		Suppose that Assumptions \ref{assumption-static} and \ref{assumption1} hold true. If $\lambda \geq \lambda_0:=\frac{\sqrt{|\mathcal{R}|}}{\tilde{\sigma}_{\min}(A)} max_{i \in \mathcal{R}} \|\nabla\E[F(\tilde{x}^*,\xi_i)]+\nabla f_0(\tilde{x}^*)\|_\infty$ where $\tilde{\sigma}_{\min}(A)$ is the minimum nonzero singular value of $A$, 
		then for the optimal solution $x^*$ of \eqref{eq6} and the optimal solution $\tilde{x}^*$ of \eqref{eq4}, we have $x^*=[\tilde{x}^*]$.
	\end{theorem}
	
	No matter how large $\lambda$ is, with a proper step size the proposed stochastic gradient method can converge to the optimal solution of \eqref{eq6} when the Byzantine agents are absent. However, the Byzantine agents bring disturbance to the optimization process, and their influence is illustrated in the following theorem.
	
	\begin{theorem} \label{theorem2}
		Suppose that Assumptions \ref{assumption-static}, \ref{assumption1}, \ref{assumption2}, \ref{assumption3} hold true. Set the step size of the proposed method given by \eqref{eq:subgradient} as $\alpha^{k}=\min\{\underline{\alpha},\frac{\overline{\alpha}}{k+1}\}$, where \black{$\underline{\alpha}=\min\limits_{i\in\mathcal{R}} \{\frac{1}{4(u_i+L_i)}\} $}, and $\overline{\alpha}>\frac{1}{\eta}$ with \black{$\eta=\min\limits_{i\in\mathcal{R}}\{\frac{2u_iL_i}{u_i+L_i}\}-\epsilon>0$ for some} $\epsilon>0$. Then, there exists a smallest integer $k_0$ satisfying $\underline{\alpha} \geq \frac{\overline{\alpha}}{k_0+1}$, such that
		\begin{align}\label{th2--eq1}
			\E\|x^{k+1} - x^*\|^2 \leq (1-\eta\underline{\alpha})^k \|x^0 - x^*\|^2 + \frac{1}{\eta}(\underline{\alpha}\Delta_0+\Delta_2), \quad \forall{k}<k_0,
		\end{align}
		and
		\begin{align}\label{th2--eq2}
			\E\|x^{k+1} - x^*\|^2 \leq \frac{\Delta_1}{k+1}+\overline{\alpha}\Delta_2, \quad \forall{k} \geq k_0.
		\end{align}
		Here we \black{take expectations over all the random variables and }define constants
		\begin{align}
			\Delta_0 = \sum\limits_{i \in \mathcal{R}} \Big(32 \lambda^2|\mathcal{R}_i|^2p+4 \lambda^2|\mathcal{B}_i|^2p +2\delta_i^2\Big), ~ \Delta_1=\max \Bigg \{\frac{\overline{\alpha}^2\Delta_0}{\eta\overline{\alpha}-1},(k_0+1)\E\|x^{k_0}-x^*\|^2+\frac{\overline{\alpha}^2\Delta_0}{k_0+1} \Bigg \}, ~ \Delta_2 = \sum\limits_{i \in \mathcal{R}}\frac{\lambda^2|\mathcal{B}_i|^2p}{\epsilon}. \nonumber
		\end{align}
	\end{theorem}
	
	Theorem \ref{theorem2} asserts that the proposed Byzantine-robust decentralized stochastic optimization method can reach a neighborhood of the optimal solution $x^*$ of \eqref{eq6}. The convergence rate is sublinear and matches the rates of Byzantine-free decentralized stochastic optimization methods \cite{dpsgd,sspm,dam,mdm}. The size of neighborhood is proportional to $p$ (the dimension of model), $\lambda^2$ (squared penalty parameter), and $\sum_{i \in \mathcal{R}} |\mathcal{B}_i|^2$ that is determined by the number of Byzantine agents and the network topology. Combining Theorems \ref{theorem1} and \ref{theorem2}, we derive the main Theorem as follows.
	
	\begin{theorem}\label{theorem3}
		Under the same conditions of Theorem \ref{theorem2}, if choosing $\lambda \geq \lambda_0$, then for a sufficiently large $k \geq k_0$, we have
		\begin{align} \label{th3-1}
			\E\|x^{k+1}-[\tilde{x}^*]\|^2 \leq \frac{\Delta_1}{k+1}+\overline{\alpha}\Delta_2.	
		\end{align}
		If choosing $0<\lambda<\lambda_0$ and supposing that the difference between the optimizers of \eqref{eq6} and \eqref{eq4} is bounded by $\|x^*-[\tilde{x}^*]\|^2\leq\Delta_3$, then for a sufficiently large $k \geq k_0$, we have
		\begin{align}\label{th3-2}
			\E\|x^{k+1}-[\tilde{x}^*]\|^2 \leq \frac{2\Delta_1}{k+1}+2\overline{\alpha}\Delta_2+ 2\Delta_3.
		\end{align}
	\end{theorem}
	
	When $\lambda$ is large enough, according to Theorem \ref{theorem1}, \eqref{eq6} is equivalent to \eqref{eq4}. Therefore, the gap between $x^k$ and $x^*$ directly translates to the gap between $x^k$ and $[\tilde{x}^*]$ as in \eqref{th3-1}. However, if $\lambda$ is too large, the gap will also be large because $\Delta_2$ is proportional to $\lambda^2$. When $\lambda$ is too small, \eqref{eq6} cannot guarantee to have a consensual solution. In this case, the gap between $[\tilde{x}^*]$ and $x^*$ is unclear, but we assume that it is bounded by $\Delta_3$. Therefore, we are also able to characterize the gap between $x^k$ and $[\tilde{x}^*]$ as in \eqref{th3-2}.
	
	\subsection{Performance Analysis over Time-Varying Network}

	Analogous to Theorem \ref{theorem1}, when the network is time-varying, we can also show that as long as $\lambda$ is large enough, the penalized problem \eqref{tvn-eq2} has the same optimal solution as the Byzantine-free problem \eqref{eq4}.
	\black{Also, we define an incidence matrix to determine how large $\lambda$ should be to guarantee consensus, yet in the time-varying case it depends on the connectivity of the average network. In specific, the weighted incidence matrix $\bar{A} \in \mathbb{R}^{|\mathcal{R}| \times |\bar{\mathcal{E}}_{R}|}$ of the average network is constructed according to the $\bar{a}_e$'s. For an edge $e=(i,j) \in \bar{\mathcal{E}}_{R}$ with $i<j$, the $(i,e)$-th entry of $\bar{A}$ is $\bar{a}_e$ while the $(j,e)$-th entry of $\bar{A}$ is $-\bar{a}_e$.} 
	\begin{theorem} \label{theorem4}
		Suppose that Assumption \ref{assumption-timevary} and \ref{assumption1} hold true. If $\lambda \geq \lambda_0:=\frac{\sqrt{|\mathcal{R}|}}{\tilde{\sigma}_{\min}(\bar{A})} max_{i \in \mathcal{R}} \|\nabla\E[F(\tilde{x}^*,\xi_i)]+\nabla f_0(\tilde{x}^*)\|_\infty$ where $\tilde{\sigma}_{\min}(\bar{A})$ is the minimum nonzero singular value of the weighted node-edge incidence matrix $\bar{A}$ of $(\mathcal{R},\bar{\mathcal{E}}_{R})$, then for the optimal solution $x^*$ of  \eqref{tvn-eq2} and the optimal solution $\tilde{x}^*$ of \eqref{eq4}, \black{almost surely } we have $x^* = [\tilde{x}^*]$.
	\end{theorem}
	
	Different to the static case where the critical parameter $\lambda_0$ depends on the node-edge incidence matrix $A$ of $(\mathcal{R}, \mathcal{\mathcal{E}_{R}})$, here $\lambda_0$ depends on the weighted node-edge incidence matrix $\bar{A}$ of $(\mathcal{R}, \mathcal{\bar{\mathcal{E}}_{R}})$. When the time-varying network degenerates to the static one, we have $\bar{A} = A$ such that Theorem \ref{theorem4} coincides with Theorem \ref{theorem1}.
	
	\black{In the time-varying network model, define $\mathcal{R}_i$ as the set of all possible neighbors of agent $i$. To be specific, $\mathcal{R}_i = \{j: j\in \mathcal{R}_i(\zeta^k)\text{ for some }\zeta^k \text{ or }j\in\mathcal{R}_i(\bar\zeta)\text{ for some }\bar\zeta\}$.
		With $\mathcal{R}_i$, we are ready to establish the following theorem which implies} the convergence property of the proposed stochastic subgradient method under Byzantine attacks.
	\begin{theorem} \label{theorem5}
		Suppose that Assumptions \ref{assumption-timevary}, \ref{assumption1}, \ref{assumption2}, \ref{assumption3} hold true. Set the step size of the proposed method given by \eqref{tvn-eq6} as $\alpha^{k}=\min\{\underline{\alpha},\frac{\overline{\alpha}}{k+1}\}$, where \black{$\underline{\alpha}=\min_{i\in\mathcal{R}} \{\frac{1}{4(u_i+L_i)}\} $, and $\overline{\alpha}>\frac{1}{\eta}$ with $\eta=\min\limits_{i\in\mathcal{R}}\{\frac{2u_iL_i}{u_i+L_i}\}-\epsilon>0$ for some $\epsilon>0$.} 
		Then, there exists a smallest integer $k_0$ satisfying $\underline{\alpha} \geq \frac{\overline{\alpha}}{k_0+1}$, such that \black{almost surely for any Byzantine attack, we have}
		\begin{align}\label{th6--eq1}
			\E\|x^{k+1} - x^*\|^2 \leq (1-\eta\underline{\alpha})^k \|x^0 - x^*\|^2 + \frac{1}{\eta}(\underline{\alpha}\Delta_4+\Delta_6), \quad \forall{k}<k_0,
		\end{align}
		and
		\begin{align}\label{th6--eq2}
			\E\|x^{k+1} - x^*\|^2 \leq \frac{\Delta_5}{k+1}+\overline{\alpha}\Delta_6, \quad \forall{k} \geq k_0.
		\end{align}
		Here \black{the expectation $\mathbb{E}$ is taken over all the random samples $\xi^k_i\sim \mathcal{D}_i$ and the constants are defined as}
		\begin{align}
			\Delta_4 = \sum\limits_{i \in \mathcal{R}} \Big(\black{32} \lambda^2|\mathcal{R}_i|^2p + 4 \lambda^2|\mathcal{B}_i|^2p +2\delta_i^2 \Big), \quad \Delta_5 =\max \Bigg \{\frac{\overline{\alpha}^2\Delta_4}{\eta\overline{\alpha}-1},(k_0+1)\E\|x^{k_0}-x^*\|^2+\frac{\overline{\alpha}^2\Delta_4}{k_0+1} \Bigg \}, \nonumber
		\end{align}
		and
		\begin{align}
			\Delta_6= \sum\limits_{i \in \mathcal{R}}\bigg(\frac{2\lambda^2|\mathcal{B}_i|^2p}{\epsilon} + \frac{8\lambda^2 |\mathcal{R}_i|^2 p}{\epsilon}\bigg) \nonumber.
		\end{align}
	\end{theorem}
	
	With Theorem \ref{theorem5}, we conclude that the proposed Byzantine-robust decentralized stochastic optimization method can also reach a neighborhood of the optimal solution $x^*$ of \eqref{tvn-eq2} when the underlying network is time-varying. Combining Theorems \ref{theorem4} and \ref{theorem5}, we get the main theorem as follows.
	\begin{theorem}\label{theorem6}
		Under the same conditions of Theorem \ref{theorem5}, if choosing $\lambda \geq \lambda_0$, then for a sufficiently large $k \geq k_0$, \black{almost surely} we have
		\begin{align} \label{th6-1}
			\E\|x^{k+1}-[\tilde{x}^*]\|^2 \leq \frac{\Delta_5}{k+1}+\overline{\alpha}\Delta_6.	
		\end{align}
		If choosing $0<\lambda<\lambda_0$ and supposing that the difference between the optimizers of \eqref{tvn-eq2} and \eqref{eq4} is bounded by $\|x^*-[\tilde{x}^*]\|^2\leq\Delta_7$, then for a sufficiently large $k \geq k_0$, \black{almost surely} we have
		\begin{align}\label{th6-2}
			\E\|x^{k+1}-[\tilde{x}^*]\|^2 \leq \frac{2\Delta_5}{k+1}+2\overline{\alpha}\Delta_6+ 2\Delta_7.
		\end{align}
	\end{theorem}
	%
	
	\section{Numerical Experiments}
	\label{sec:4}
	In this section, we conduct several numerical experiments to demonstrate the robustness of our proposed method to Byzantine attacks over static and time-varying networks.
	
	\subsection{Numerical Experiments over Static Network}
	\label{sec:4a}
	
	Consider a static Erdos-Renyi network consisting of $n=30$ agents, where every edge $e$ is activated with probability $p_e = 0.7$. Randomly choose $b$ agents to be Byzantine, but guarantee that the network of regular agents is connected.
	
	\black{We consider the softmax regression problem, which is defined as
		\begin{align}
			F(\tilde{x}) = - \frac{1}{N} \sum_{l=1}^{N} \sum_{m=0}^{M-1} \Big(I(w^{(l)} = m) \ln \Big(\frac{exp((\tilde{x})_{m}^T v^{(l)})}{\sum_{j=0}^{M-1} exp((\tilde{x})_j^T v^{(l)})}\Big)\Big) \nonumber.
		\end{align}
		Here $N$ and $M$ are the numbers of data samples and data classes, respectively. Further, $(v^{(l)}, w^{(l)})$ is the $l$-th data sample with $v^{(l)} \in \R^{p/M}$ and $w^{(l)} \in \R$, and $I(w^{(l)} = m)$ is the indicator function with $I(w^{(l)} = m) = 1$ if $w^{(l)} = m$ and $I(w^{(l)} = m) = 0$ otherwise. The model parameter is $\tilde{x} \in \R^p$ and $(\tilde{x})_j \in \R^{p/M}$ is the $j$-th block of $\tilde{x}$. The dataset is MNIST, which contains $M = 10$ handwritten digits from 0 to 9, with 60,000 training images and 10,000 testing images whose dimensions are $p/M = 784$.} In the i.i.d. case, we randomly and evenly distribute the training images to all the agents. In the non-i.i.d. case, we let every three agents evenly split the training images of one digit. We use the softmax regression with regularization term $f_0(\tilde{x}) = \frac{0.01}{2} \|\tilde{x}\|_2^2$ to learn the model. At the testing stage, we randomly choose one regular agent and use its local model to calculate the classification accuracy. Also, we calculate the variance of regular agents' local models to quantify the level of consensus.
	
	The benchmark methods are DPSGD \cite{dpsgd}, as well as the stochastic versions of ByRDiE \cite{ByRDiE} and BRIDGE \cite{BRIDGE} (denoted by ByRDiE-S and BRIDGE-S, respectively). In DPSGD, the mixing matrix is constructed following the equal neighbor weights rule \cite{equal-weight}. \black{In ByRDiE-S, the coordinates of the model are updated sequentially. We set the number of inner-loop iterations to update every coordinate to be 1, as suggested by \cite{ByRDiE}.} For fair comparison, in ByRDiE-S one iteration refers to that all the coordinates have been updated once. Step sizes of the benchmark methods are hand-tuned to the best. \black{In all the compared methods, the batch-size is set as 32.}
	
	
	\vspace{-1em}
	
	\begin{figure}[H]
		\centering
		\includegraphics[scale=0.4]{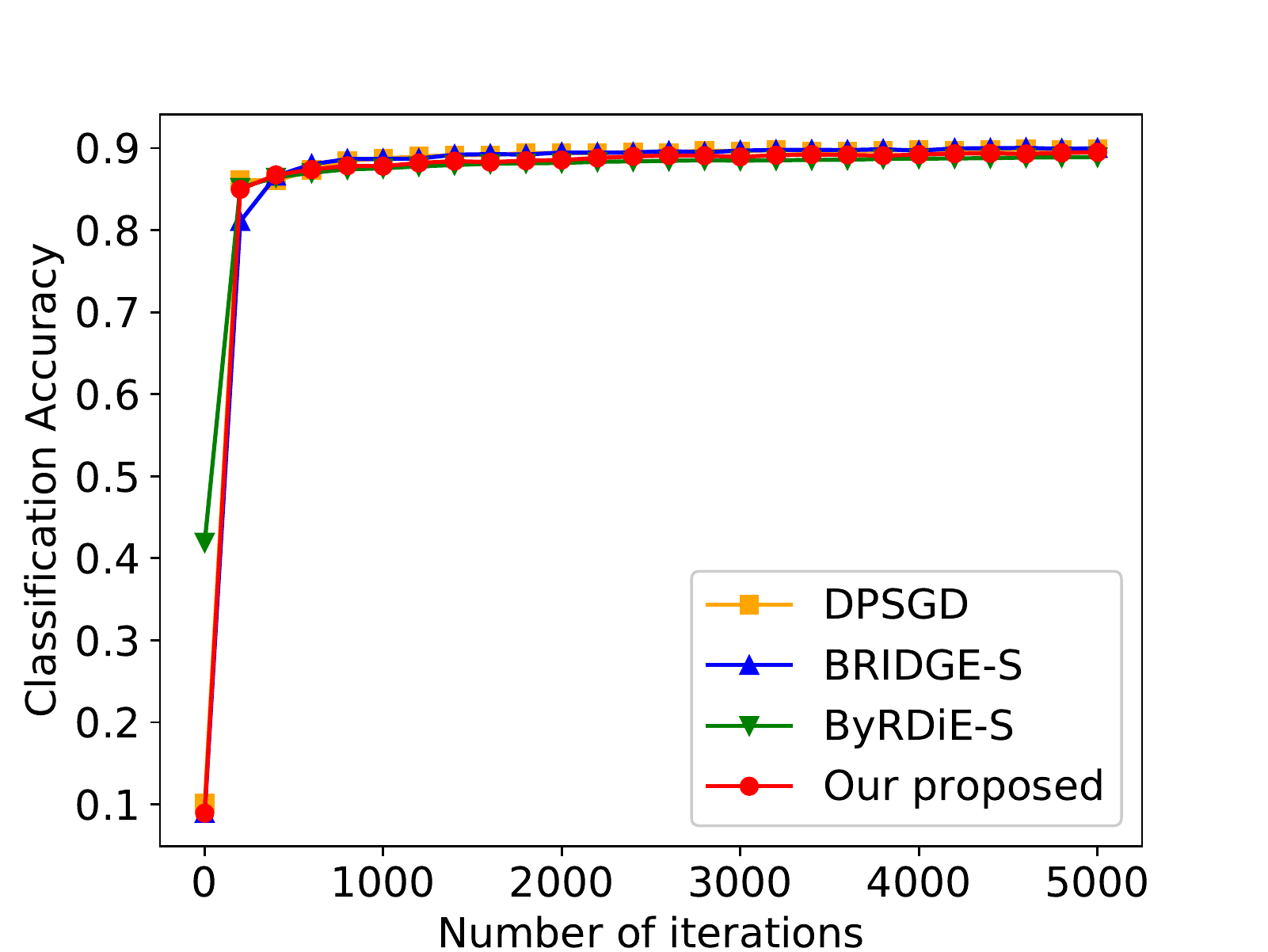} \hspace{-1.5em}
		\includegraphics[scale=0.4]{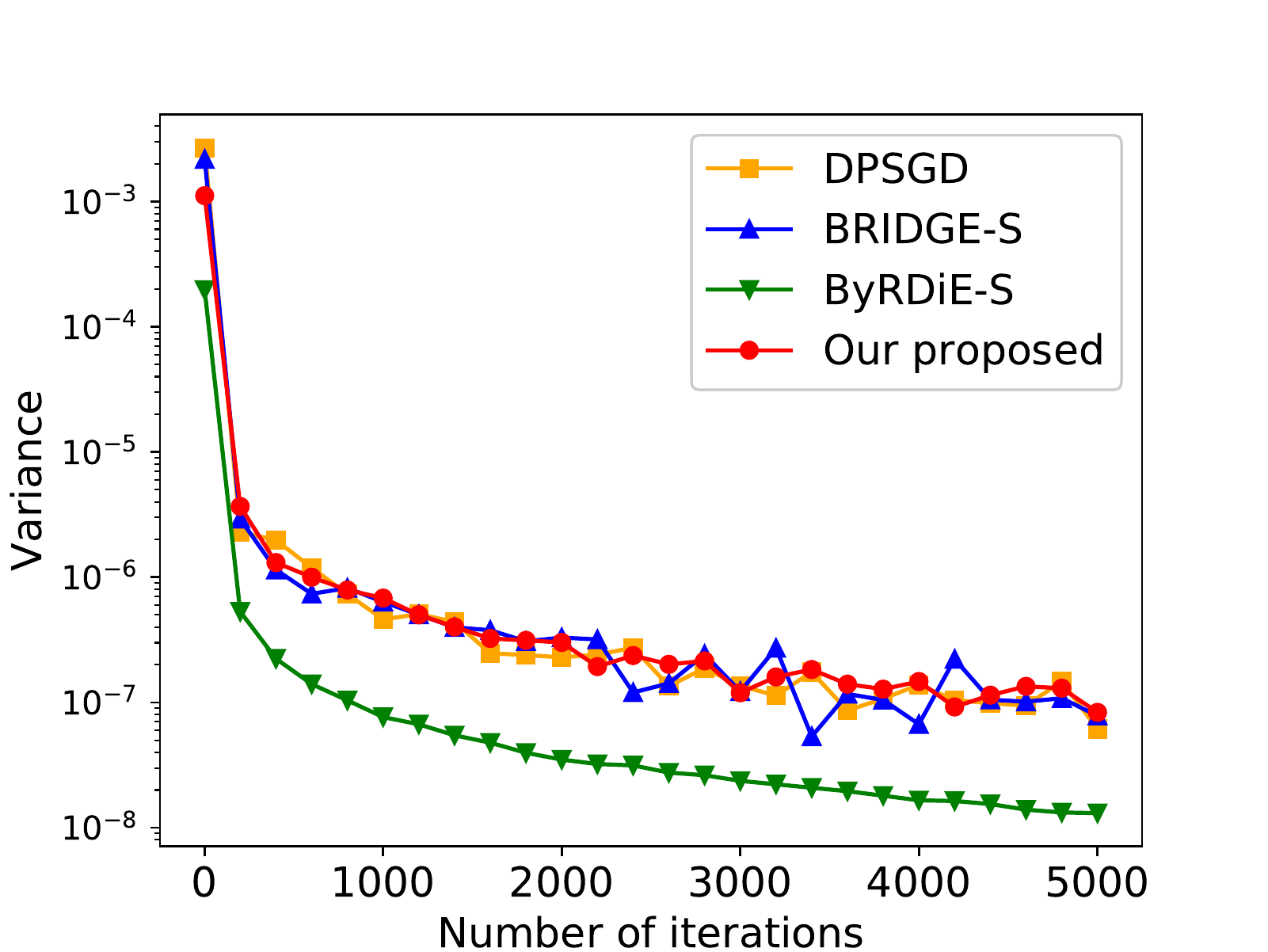}
		\caption{Classification accuracy and variance of regular agents' local models without Byzantine attacks.} \label{WA}
	\end{figure}

    \vspace{-2em}

    \begin{figure}[H]
    	\centering
    	\includegraphics[scale=0.4]{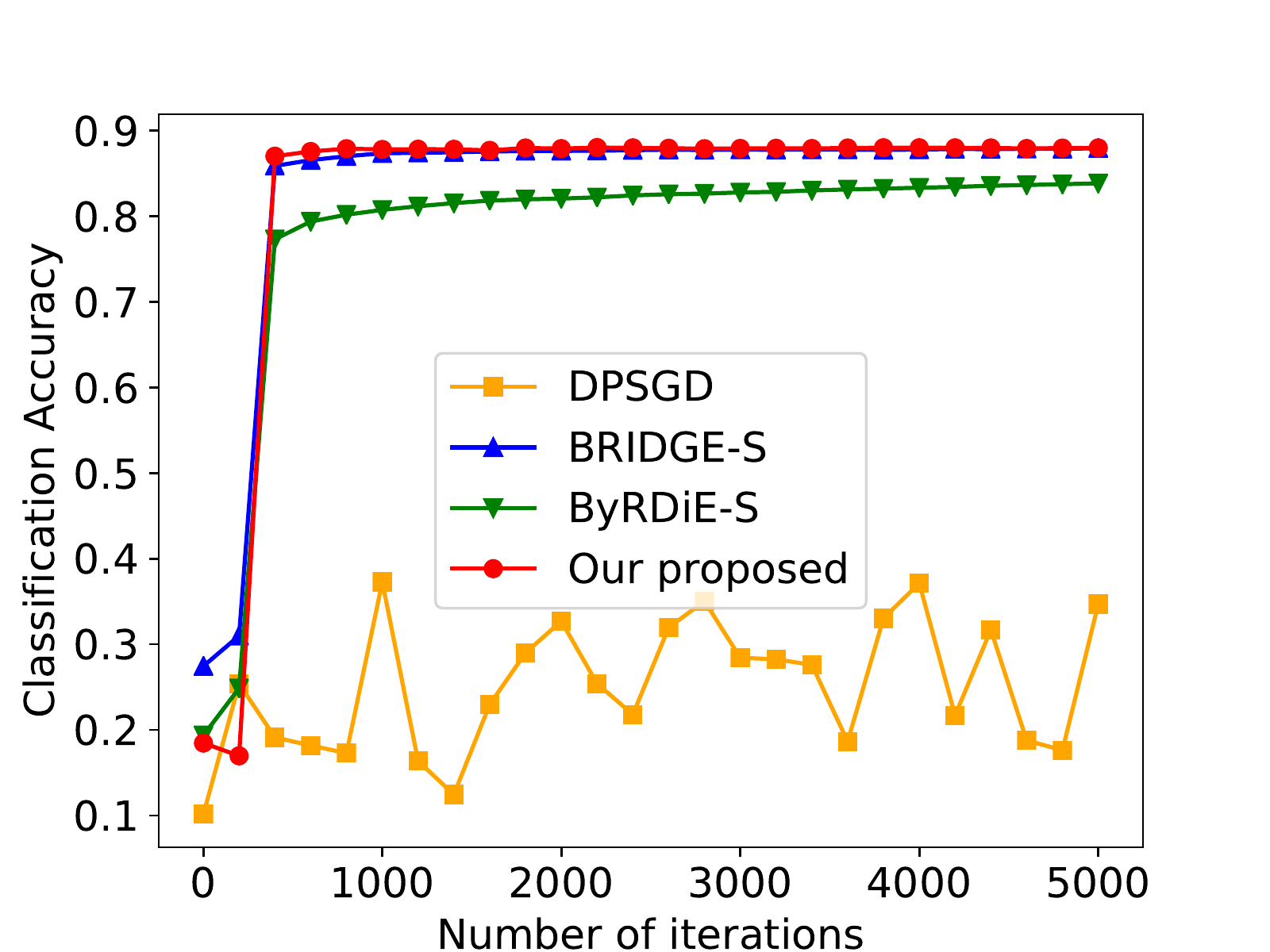} \hspace{-1.5em}
    	\includegraphics[scale=0.4]{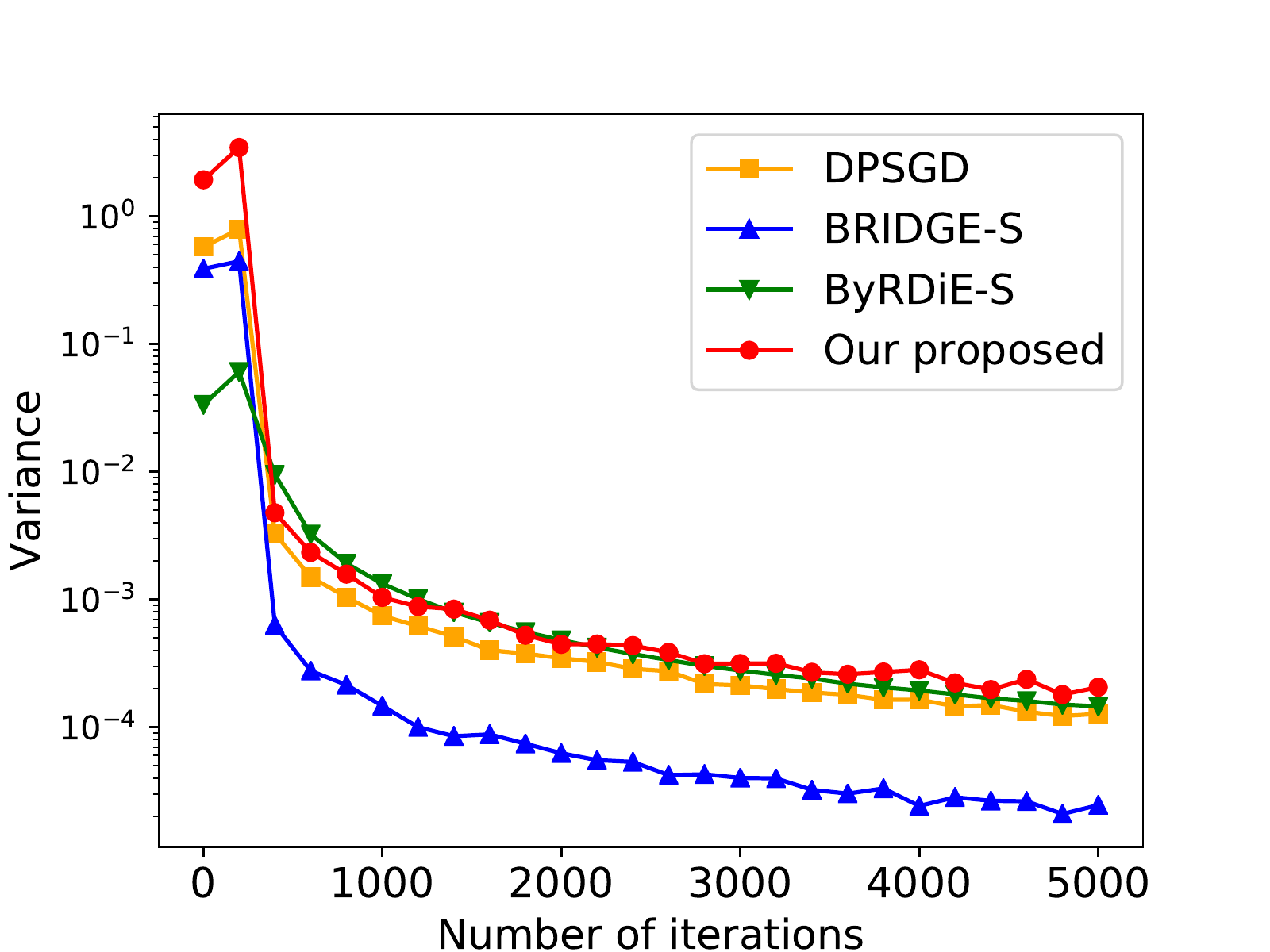}
    	\caption{Classification accuracy and variance of regular agents' local models under zero-sum attacks.} \label{ZS}
    \end{figure}
	
	\vspace{-2em}
	
	\begin{figure}[H]
		\centering
		\includegraphics[scale=0.4]{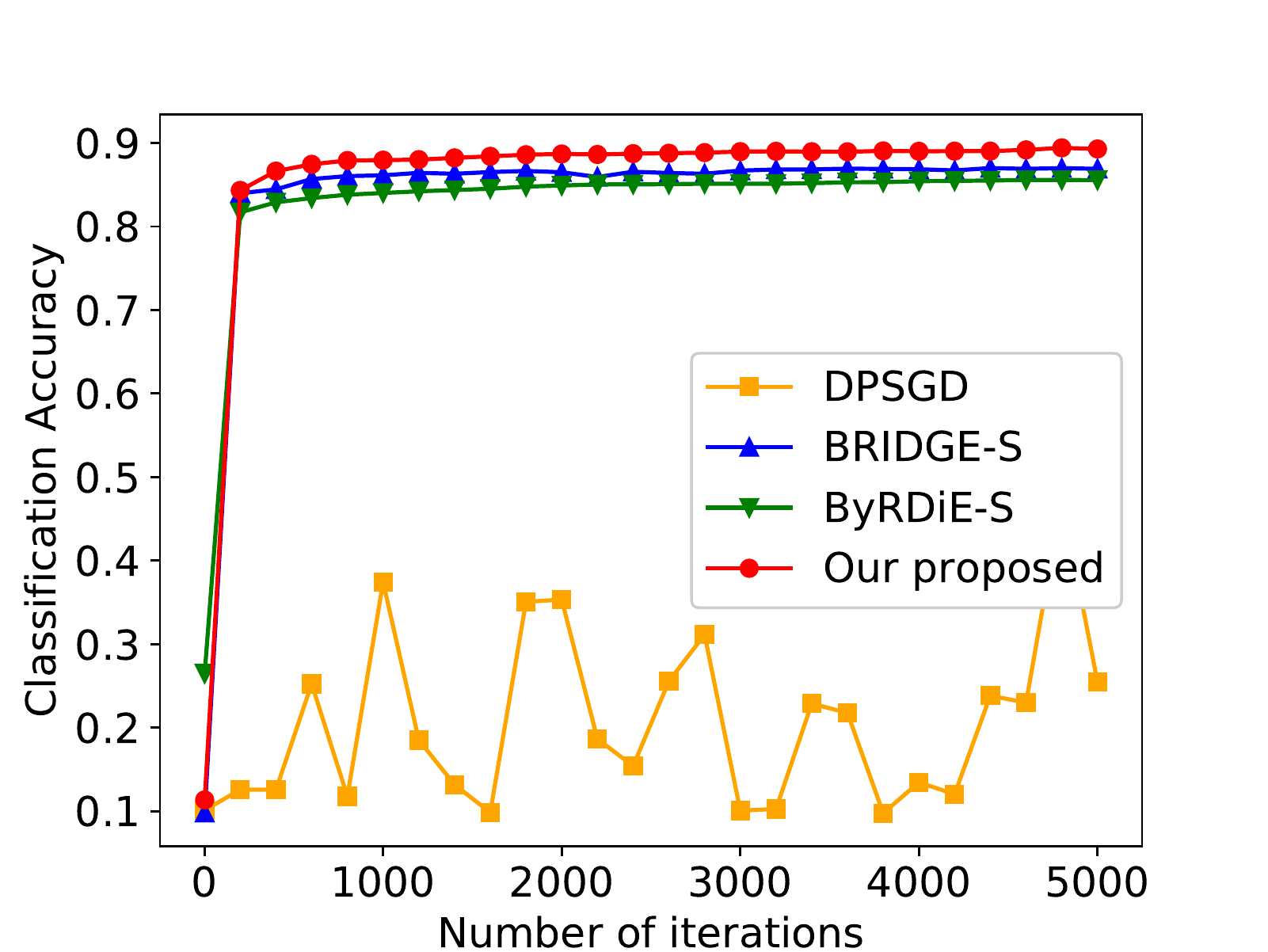} \hspace{-1.5em}
		\includegraphics[scale=0.4]{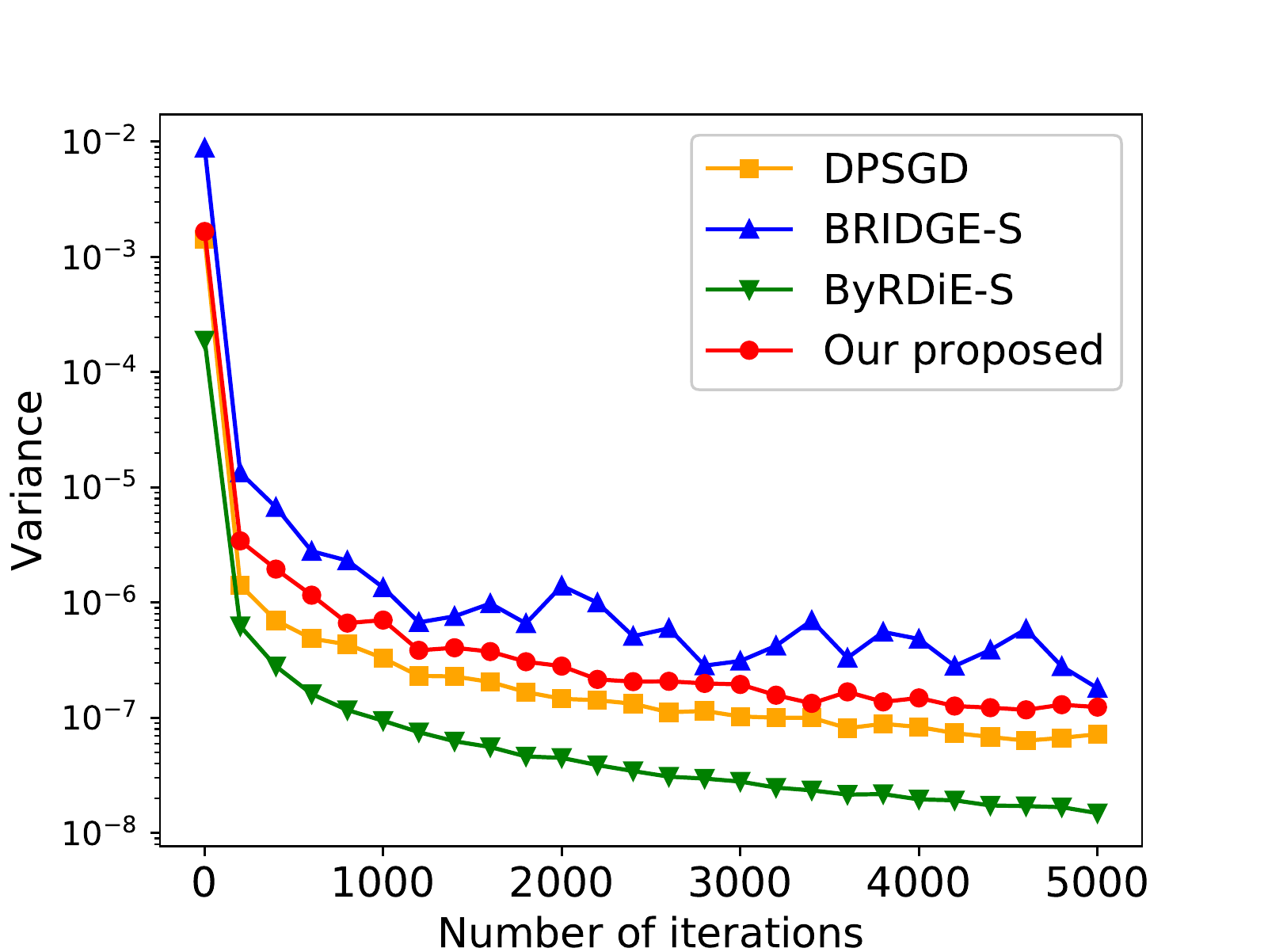}
		\caption{Classification accuracy and variance of regular agents' local models under same-value attacks.} \label{SV}
	\end{figure}
	
	\vspace{-2em}
	
	\begin{figure}[H]
		\centering
		\includegraphics[scale=0.4]{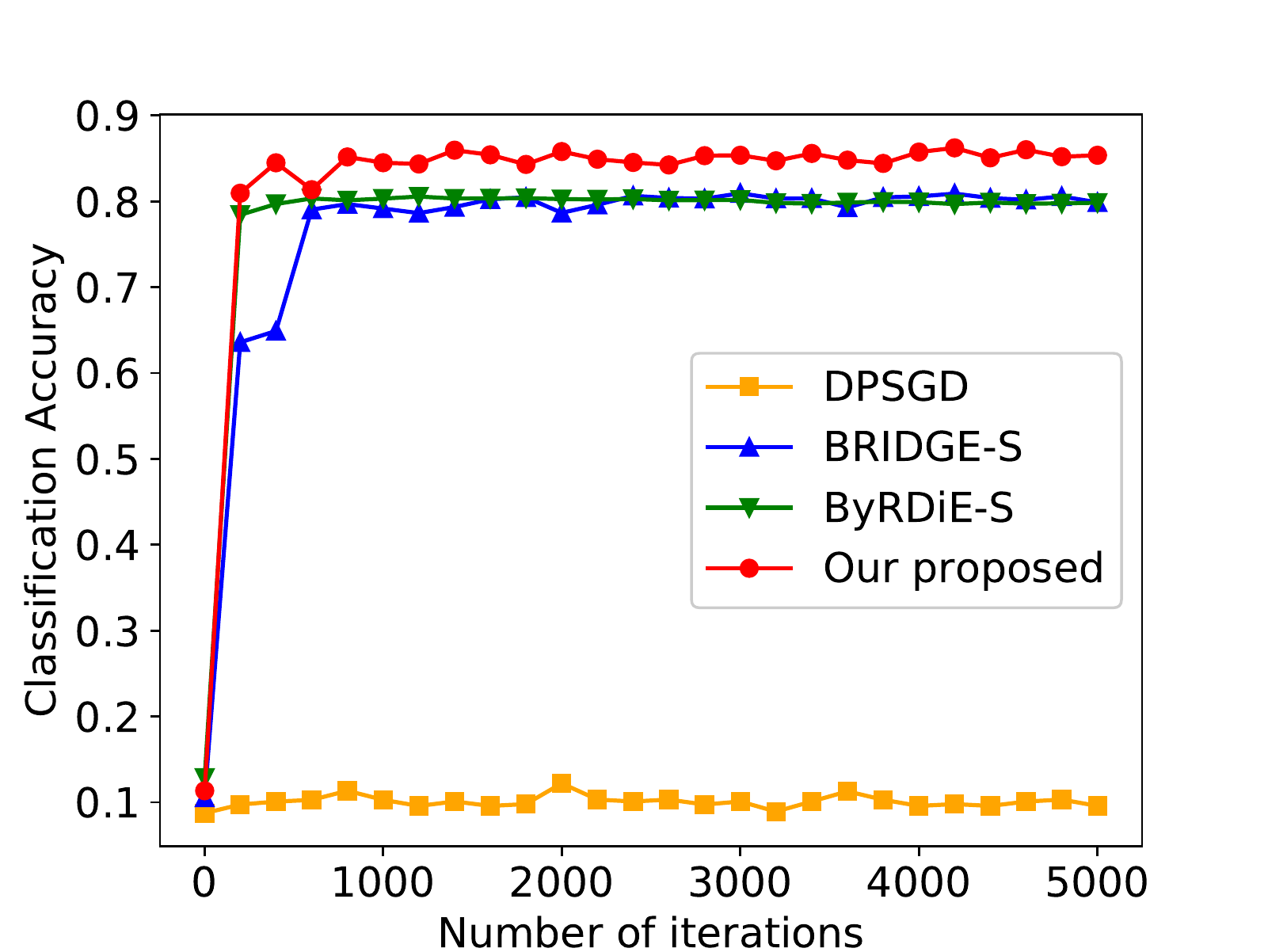} \hspace{-1.5em}
		\includegraphics[scale=0.4]{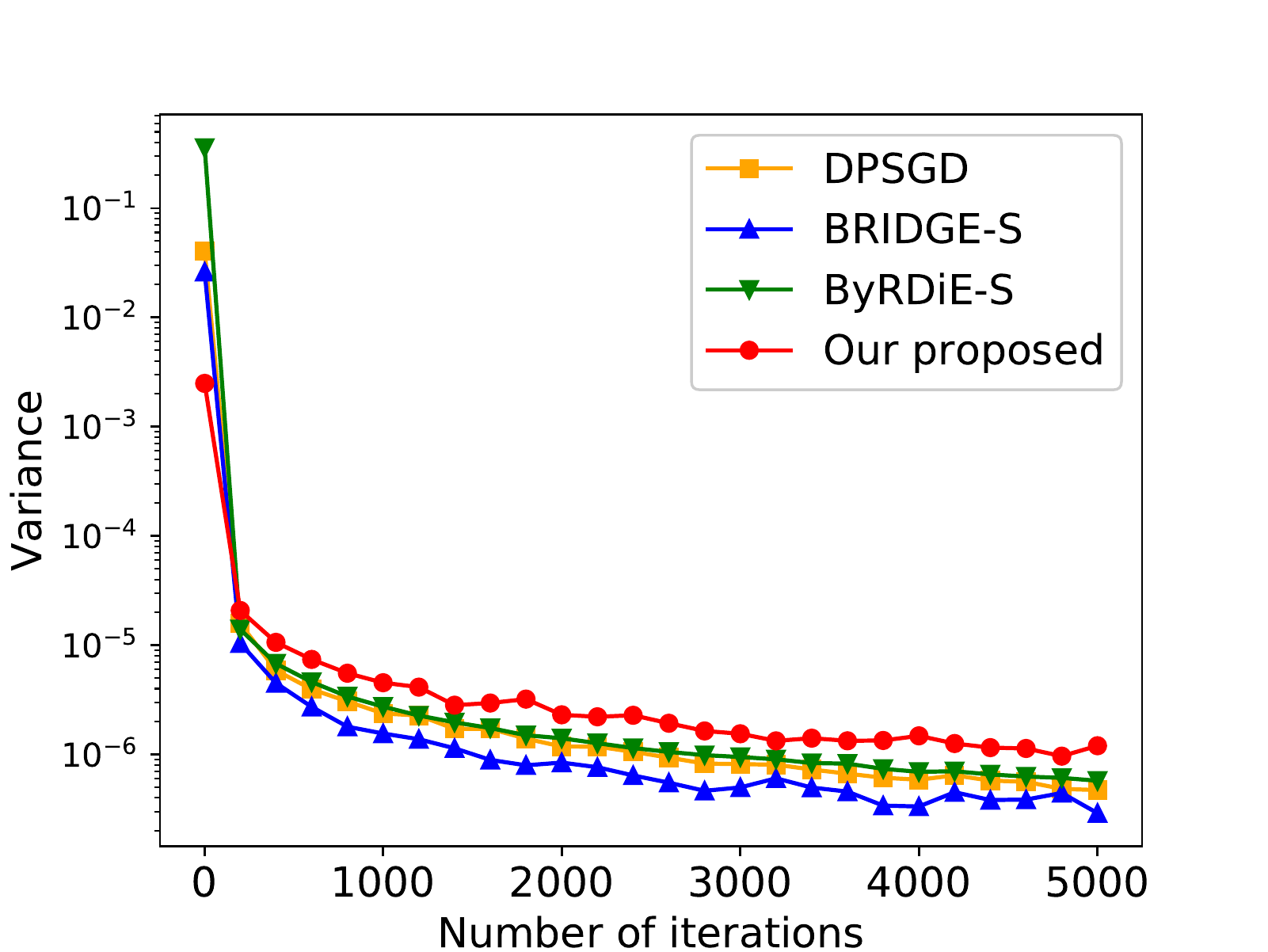}
		\caption{Classification accuracy and variance of regular agents' local models under sign-flipping attacks.} \label{SF}
	\end{figure}
	
	\noindent \textbf{Without Byzantine Attacks}. When the number of Byzantine agents is $b=0$, all the methods perform well in terms of both classification accuracy and level of consensus, as depicted in Figure \ref{WA}. In our proposed method, we set the penalty parameter as $\lambda=0.005$ and the step size as $\alpha^k=\frac{0.3}{\sqrt{k+1}}$. \black{Because of the sensitivity of the $O(\frac{1}{k})$ step size to its initial value, we use the $O(\frac{1}{\sqrt{k}})$ step size in the numerical experiments.}
	
	\noindent \textbf{Zero-sum Attacks.} Let the number of Byzantine agents be $b = 3$. Every Byzantine agent $j \in \mathcal{B}$ sends $z_j^k = - \frac{\sum_{l \in \mathcal{R}_i } x_l^k}{|\mathcal{B}_i|}$ to its regular neighbor $i \in \mathcal{R}$, such that the received messages of $i$ are summed to zero. In our proposed method, the penalty parameter is $\lambda = 0.001$ and the step size is $\alpha^k = \frac{0.9}{\sqrt{k+1}}$. With results depicted in Figure \ref{ZS}, we observe that DPSGD fails because of its vulnerability to Byzantine attacks. Our proposed method, BRIDGE-S and ByRDiE-S are robust to the zero-sum attacks.
	
	\noindent \textbf{Same-value Attacks}. Let the number of Byzantine agents be $b=3$. Every Byzantine agent $j \in \mathcal{B}$ sends $z_j^k=c\mathbf{1}$ to its neighbors. Here $\mathbf{1} \in \R^p$ is an all-one vector and $c$ is a constant which we set as 100. In our proposed method, the penalty parameter is $\lambda=0.01$ and the step size is $\alpha^k=\frac{0.28}{\sqrt{k+1}}$. As shown in Figure \ref{SV}, DPSGD fails and our proposed method is the best among all the three Byzantine-robust methods in terms of classification accuracy. Its variance is higher than that of ByRDiE-S, but small enough such that all the regular agents have high classification accuracies.
	
	\noindent \textbf{Sign-flipping Attacks}. Let the number of Byzantine agents be $b=3$. Every Byzantine agent $j \in \mathcal{B}$ first calculates its true model, and then multiplies it with a negative constant $\gamma$ and sends to its neighbors. Here we set $\gamma=-4$. In our proposed method, the penalty parameter is $\lambda=0.0022$ and the step size is $\alpha^k=\frac{0.5}{\sqrt{k+1}}$. As shown in Figure \ref{SF}, the results are consistent with those under the same-value attacks, but the performance gain of our proposed method in terms of classification accuracy is more obvious. Note that we choose a relatively small $\lambda$ such that the level of consensus is slightly worse than those of ByRDiE-S and BRIDGE-S.
	
	\vspace{-1em}
	
	\begin{figure}[H]
		\centering
		\includegraphics[scale=0.4]{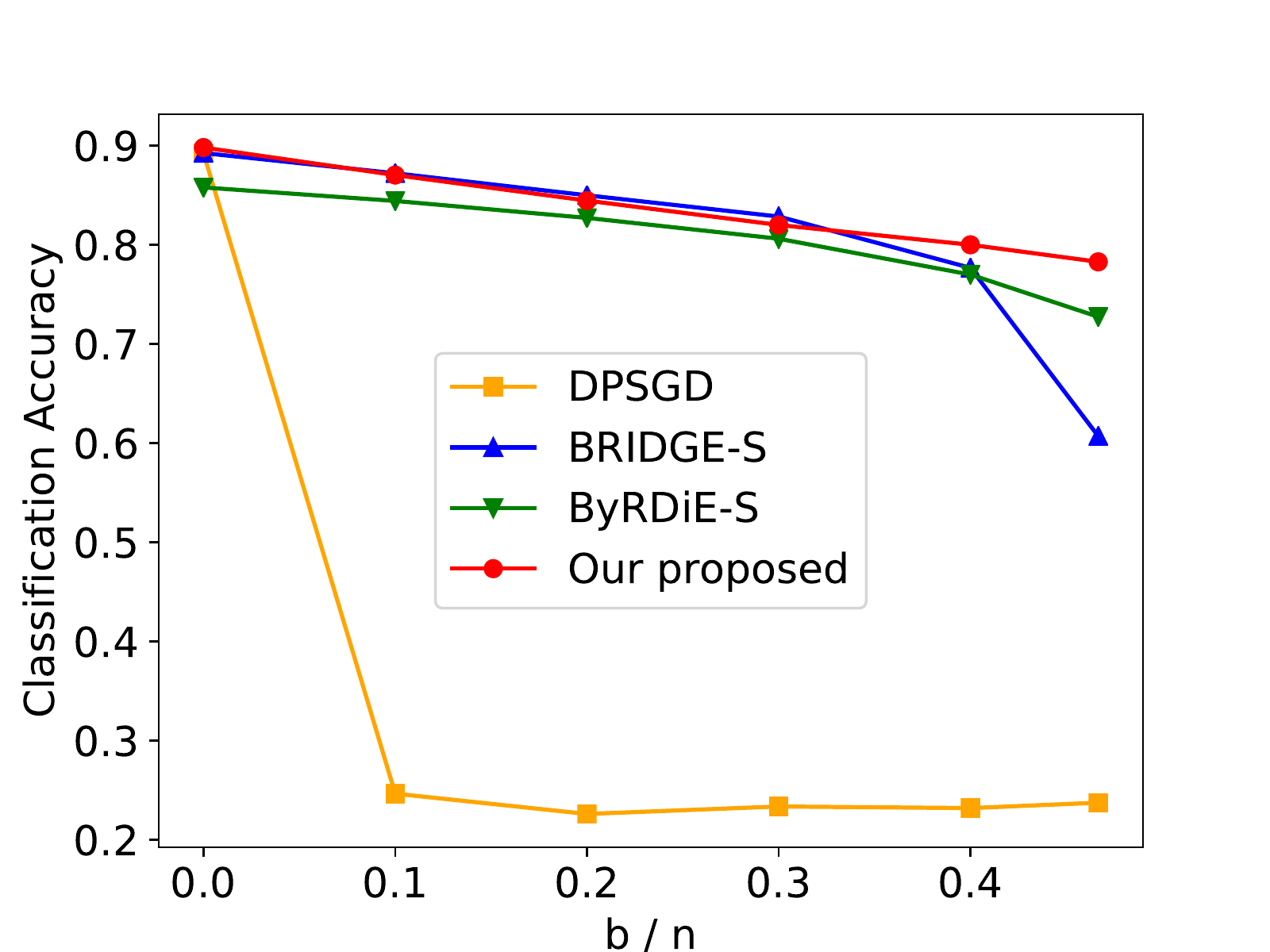} \hspace{-1.5em}
		\includegraphics[scale=0.4]{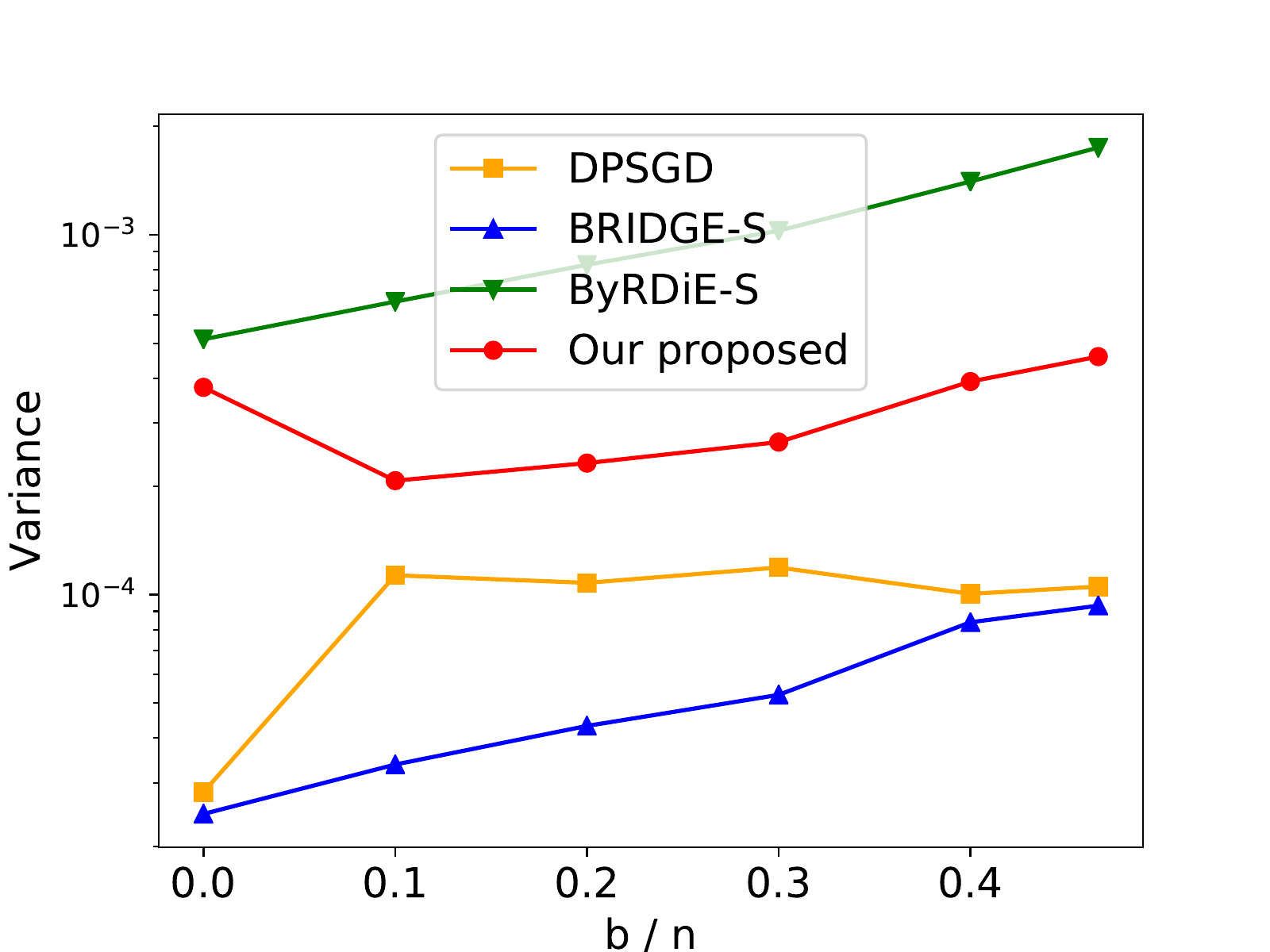}
		\caption{\black{Classification accuracy and variance of regular agents' local models with different fractions of Byzantine agents under zero-sum attacks.}} \label{Imob-iid}
	\end{figure}
	
	\vspace{-2em}
	
	\begin{figure}[H]
		\centering
		\includegraphics[scale=0.4]{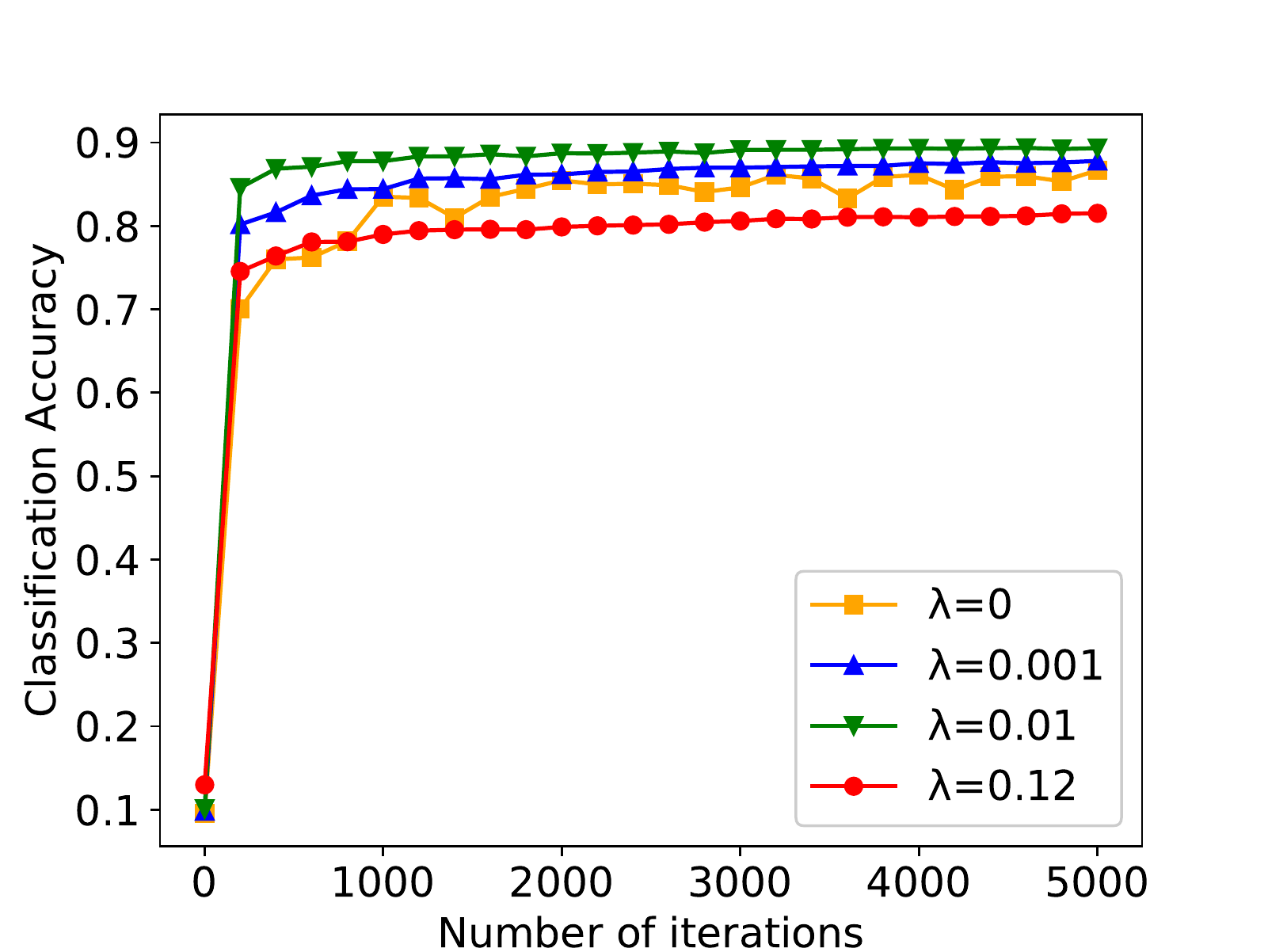} \hspace{-1.5em}
		\includegraphics[scale=0.4]{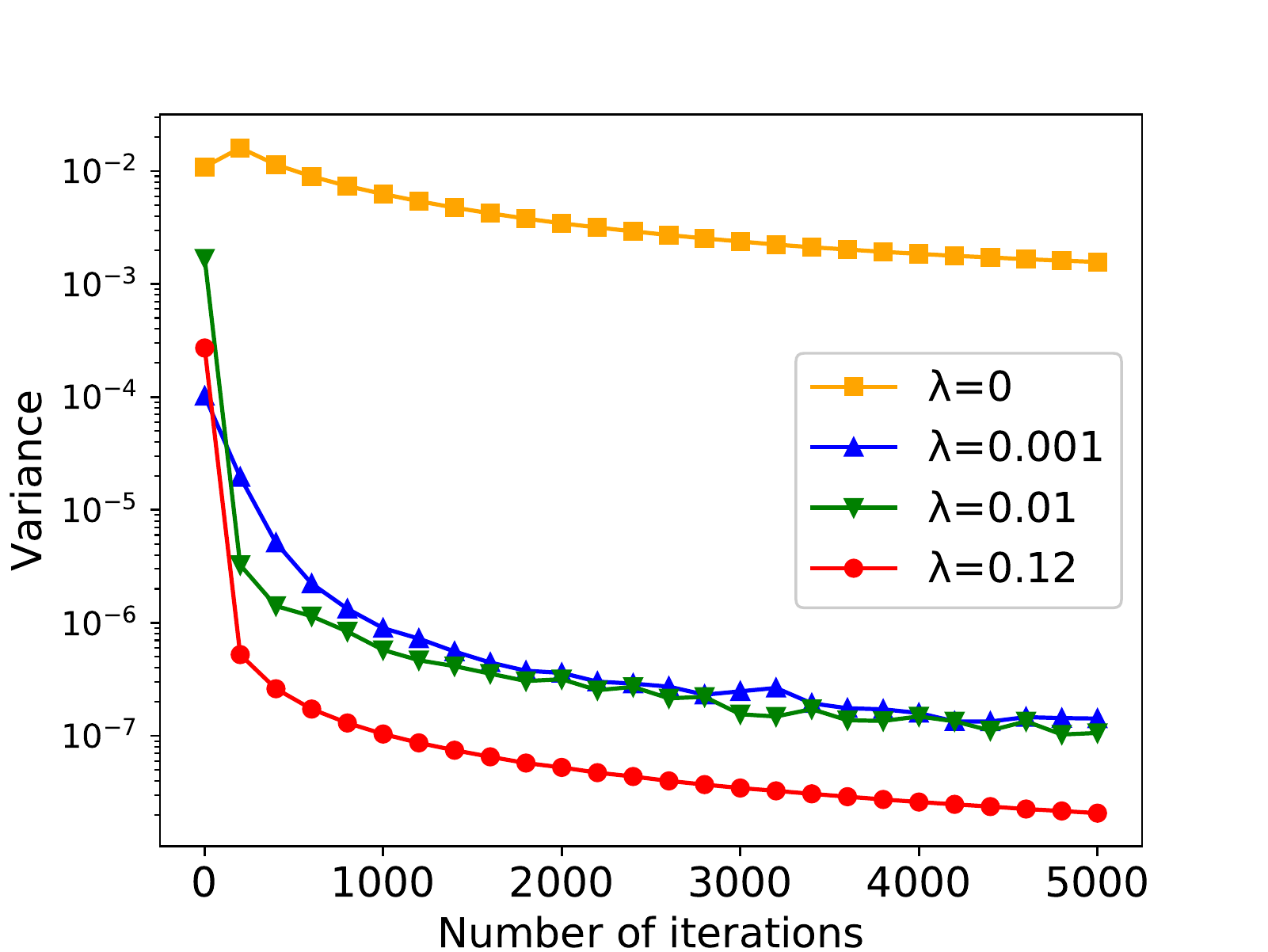}
		\caption{Classification accuracy and variance of regular agents' local models with different $\lambda$ under same-value attacks.} \label{Imopp}
	\end{figure}
	
	\vspace{-2em}
	
	\begin{figure}[H]
		\centering
		\includegraphics[scale=0.4]{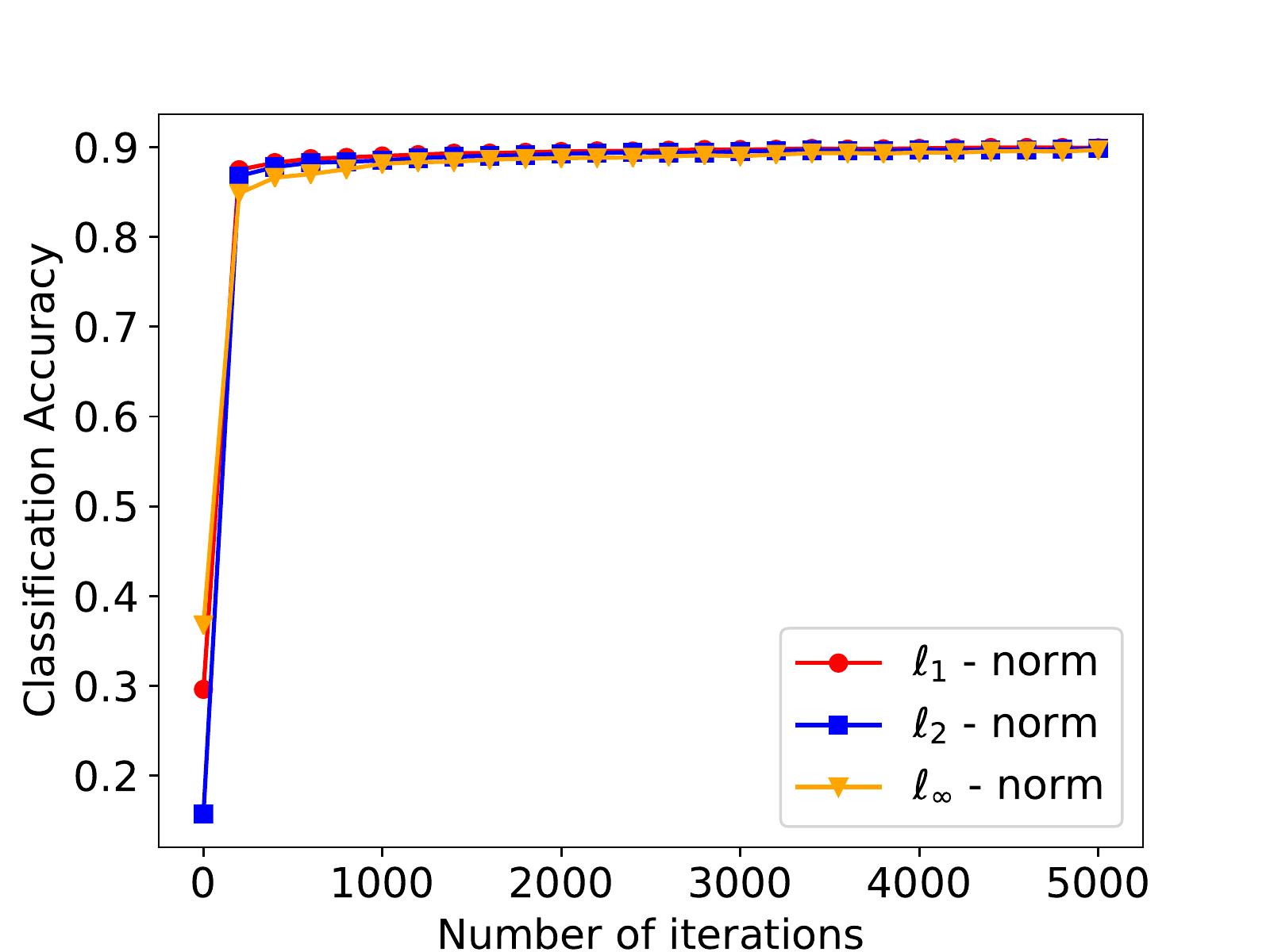} \hspace{-1.5em}
		\includegraphics[scale=0.4]{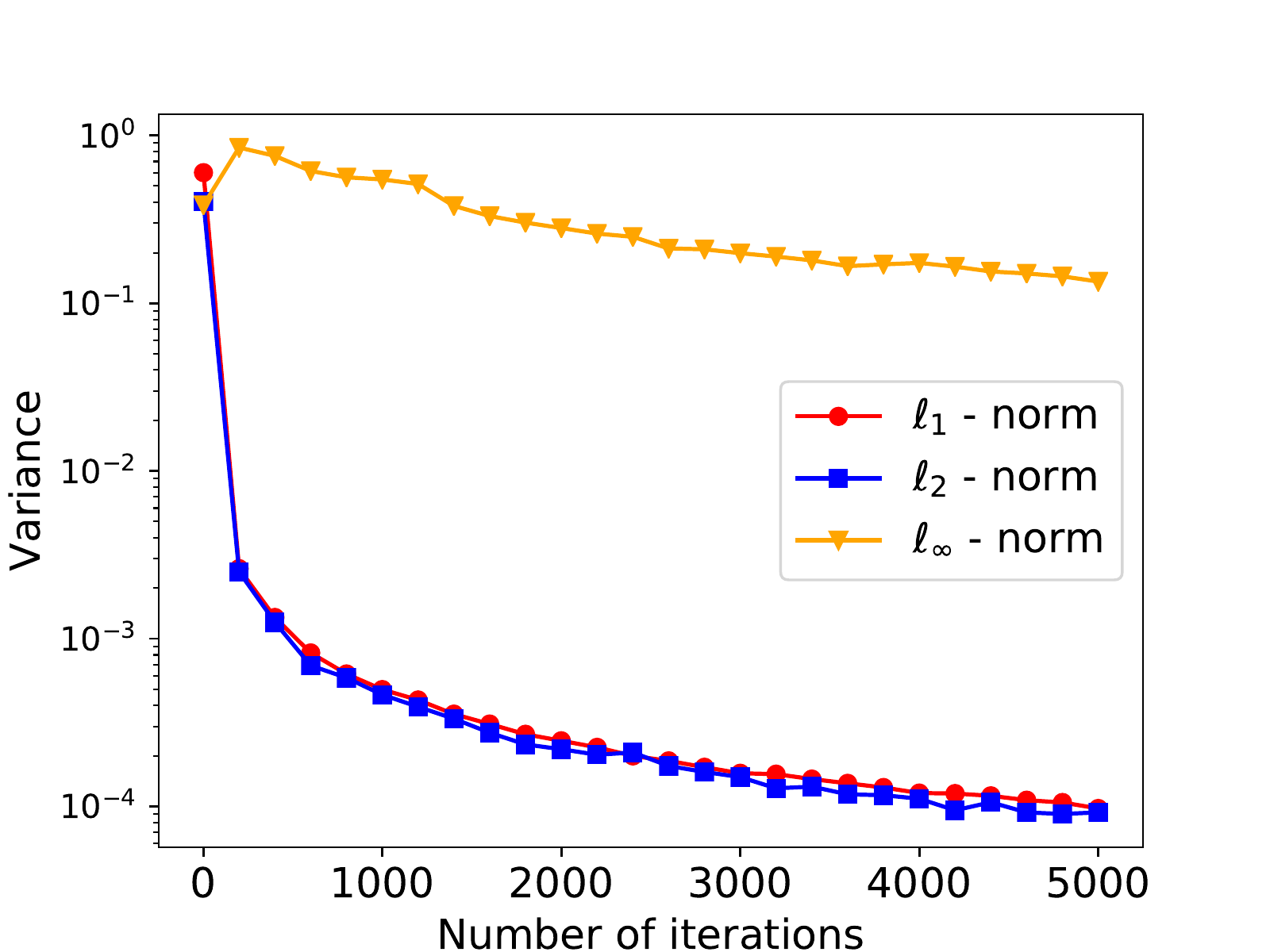}
		\caption{\black{Classification accuracy and variance of regular agents' local models based on different norms under sign-flipping attacks.}} \label{Norms}
	\end{figure}
	
	\black{
		\noindent \textbf{Impact of Fraction of Byzantine Agents}. We study the robustness of our proposed method with different fractions of Byzantine agents under the zero-sum attacks. In our proposed method, we set the penalty parameter $\lambda = 0.001$ and the step size $\alpha^k = \frac{0.9}{\sqrt{k+1}}$. Since ByRDiE-S and BRIDGE-S require $n \geq 2 b + 1$, we let no more than half of the agents be Byzantine. When the fraction of Byzantine agents increases, as depicted in Figure \ref{Imob-iid}, the performance of all the methods degrades slightly.
	}
	
	\noindent \textbf{Impact of Penalty Parameter $\lambda$}. To investigate the impact of penalty parameter $\lambda$, we choose several different values for $\lambda$ in the setting of same-value attacks with $b=3$ Byzantine agents. The step sizes are hand-tuned to the best. As shown in Figure \ref{Imopp}, larger $\lambda$ ensures better consensus, which corroborates the theoretical results in Section \ref{sec:3}. When $\lambda=0$, the level of consensus is the worst, since the agents do not communicate and learn with their own local data samples independently. However, larger $\lambda$ leads to larger gap relative to the Byzantine-free optimal solution, and hence lower classification accuracy. This observation also matches the theoretical results in Section \ref{sec:3}.
	
	\black{
		\noindent \textbf{Penalty with Different Norms}. To validate the effectiveness of TV norm penalty with different norms, we compare our proposed method with $\ell_1$, $\ell_2$  and $\ell_\infty$ norms under sign-flipping attacks. We set $\lambda=0.0022$ and $\alpha^k = \frac{0.5}{\sqrt{k+1}}$ in the $\ell_1$ case, $\lambda = 0.2$ and $\alpha^k = \frac{0.4}{\sqrt{k+1}}$ in the $\ell_2$ case, while $\lambda = 0.9$ and $\alpha^k = \frac{0.4}{\sqrt{k+1}}$ in the $\ell_\infty$ case. As shown in Figure \ref{Norms}, the accuracies of our proposed method with different norms are similar with different norms, but the result with $\ell_\infty$ norm has the biggest variance since it needs a large $\lambda$ to guarantee consensus. We deliberately choose a moderate $\lambda$ to ensure accuracy but sacrifice consensus in the $\ell_\infty$ case.
	}
	
	\vspace{-1em}
	
	\begin{figure}[H]
		\centering
		\includegraphics[scale=0.4]{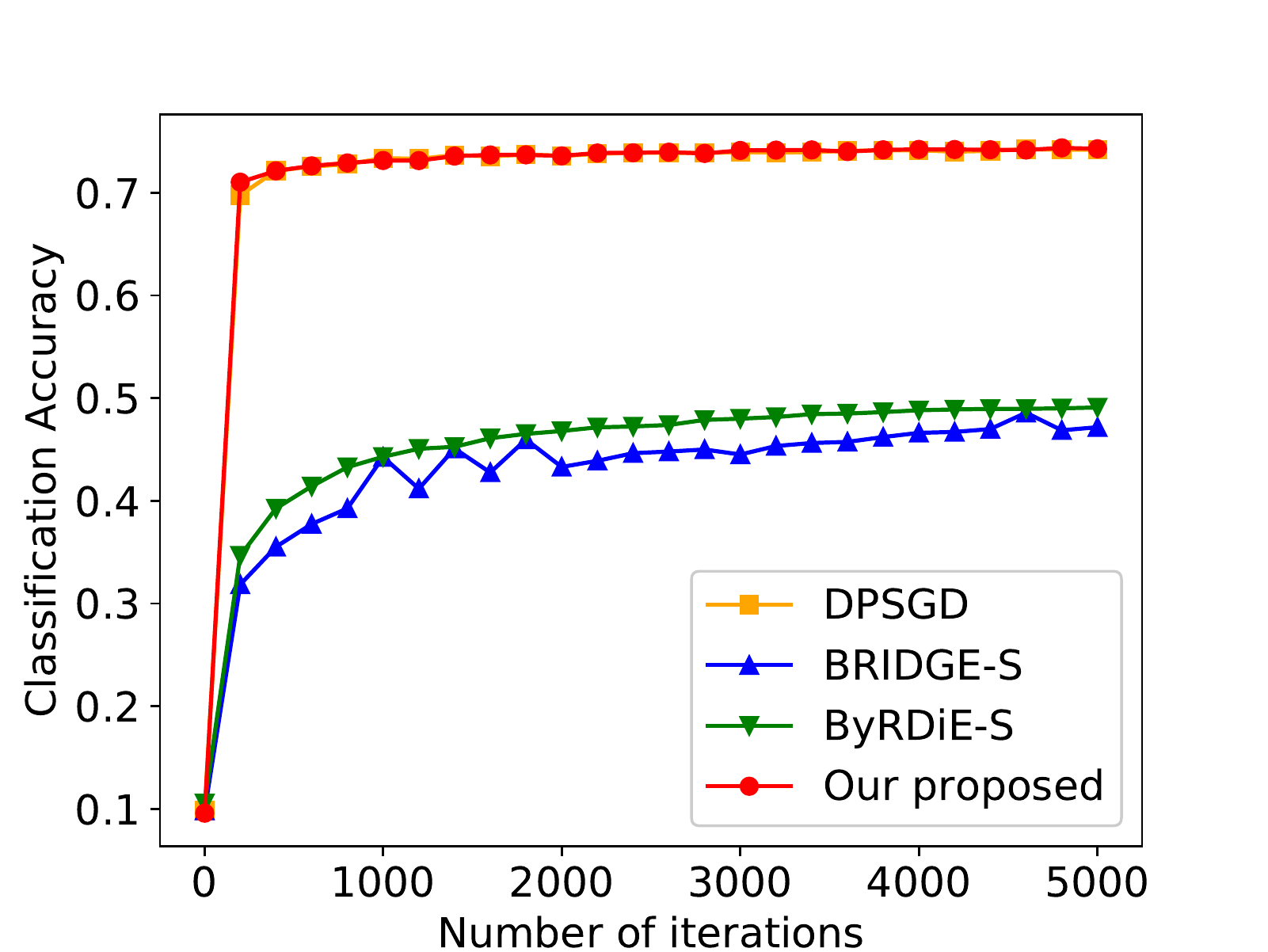} \hspace{-1.5em}
		\includegraphics[scale=0.4]{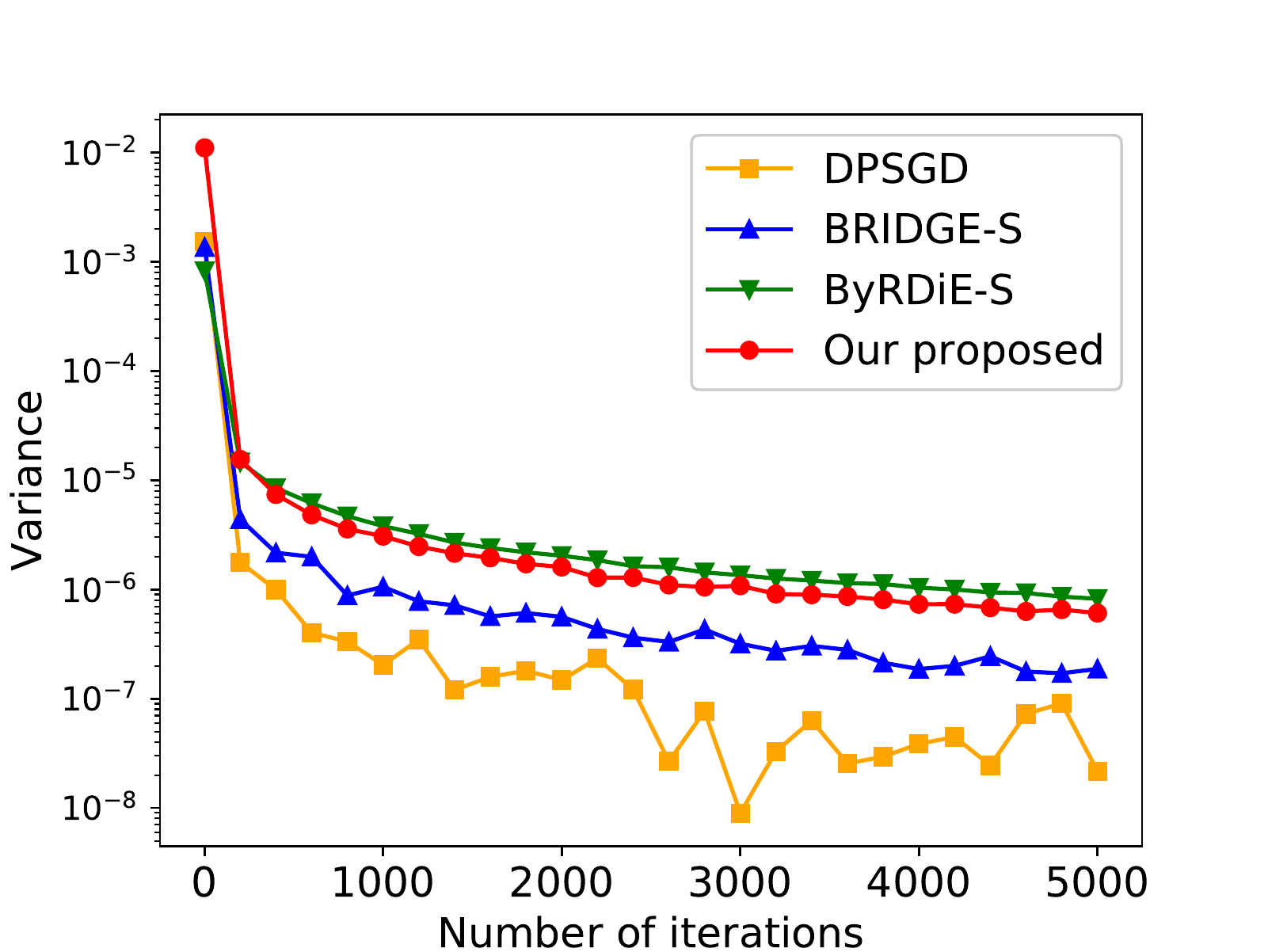}
		\caption{Classification accuracy and variance of regular agents' local models with non-i.i.d. data.} \label{Niid}
	\end{figure}
	
	\vspace{-2em}
	
	\begin{figure}[H]
		\centering
		\includegraphics[scale=0.4]{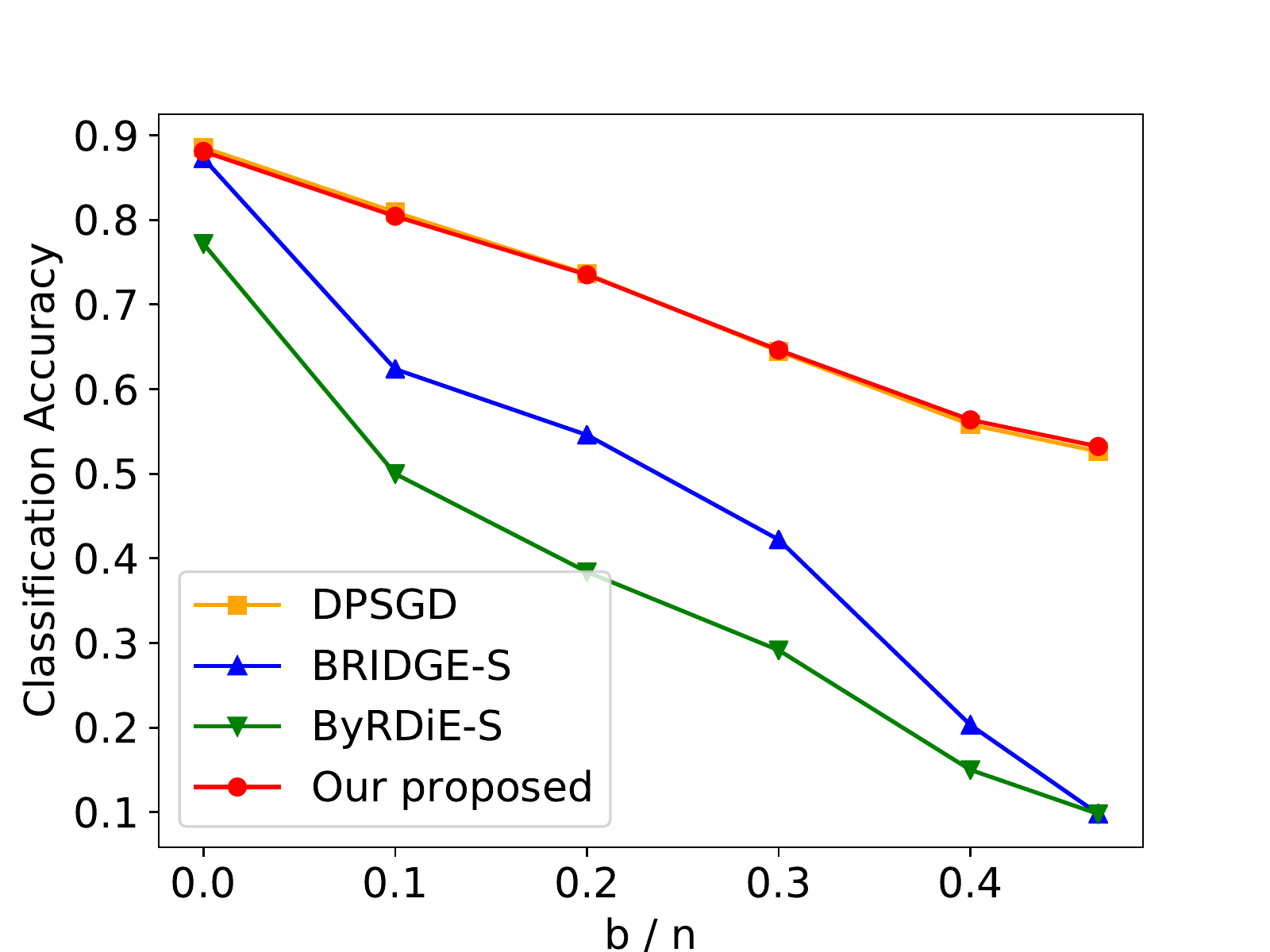} \hspace{-1.5em}
		\includegraphics[scale=0.4]{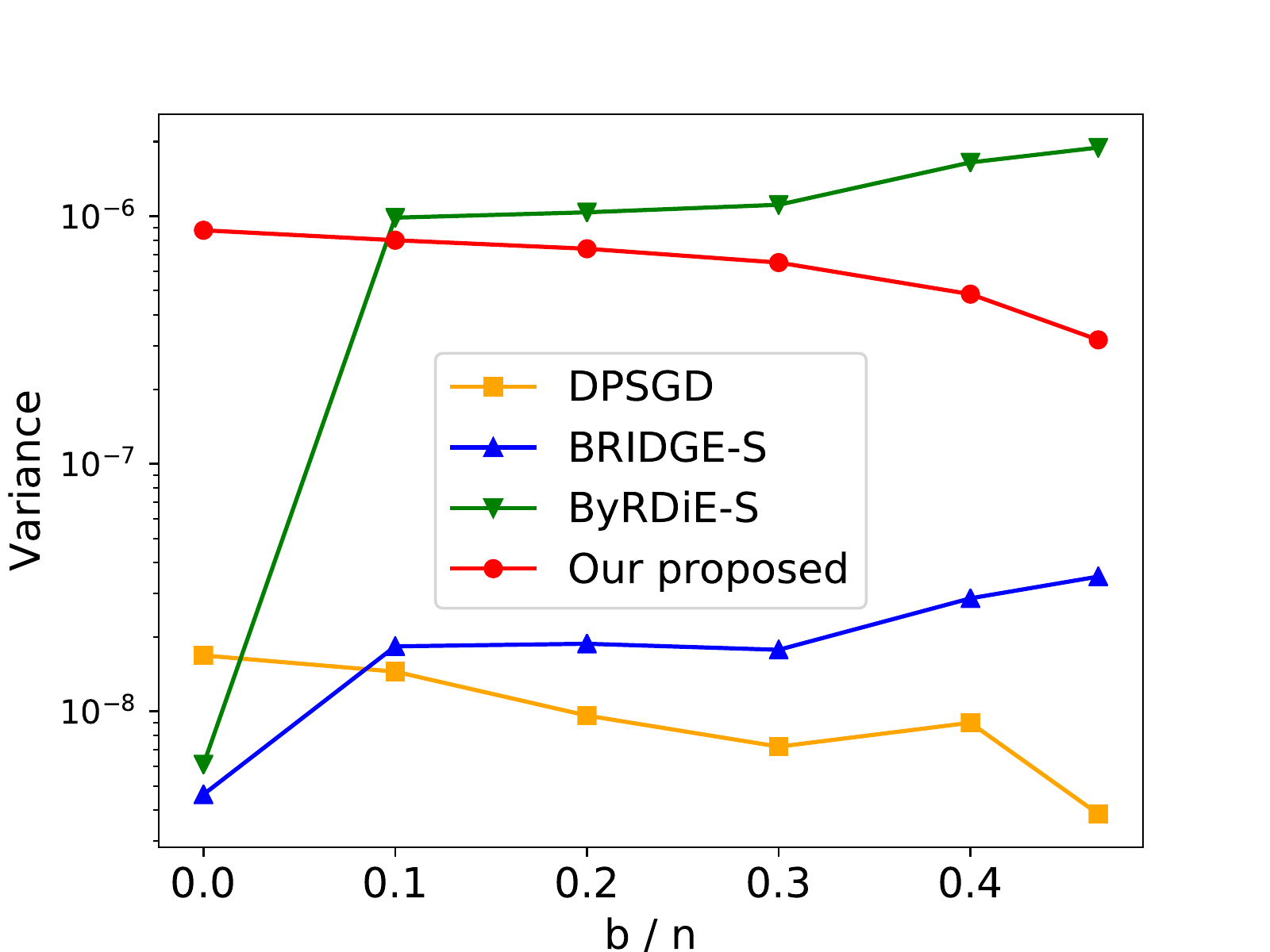}
		\caption{\black{Classification accuracy and variance of regular agents' local models with different fractions of Byzantine agents and  non-i.i.d. data.}} \label{Imob-noniid}
	\end{figure}
	
	\noindent \textbf{Non-i.i.d. Data}. Let the number of Byzantine agents be $b=6$. All the Byzantine agents copy the values of one randomly chosen regular agent, and send to their neighbors. Recall that every three agents evenly split the training images of one digit and here we deliberately let the Byzantine agents share the training images of digits 8 and 9. Therefore, information from digits 8 and 9 is totally lost and the best classification accuracy we can reach is no more than 0.8. Note that under these particularly designed attacks, DPSGD is able to reach a satisfactory classification accuracy. In our proposed method, the penalty parameter is $\lambda=0.02$ and the step size is $\alpha^k=\frac{0.4}{\sqrt{k+1}}$. As shown in Figure \ref{Niid}, our proposed method almost coincides with DPSGD with respect to classification accuracy. ByRDiE-S and BRIDGE-S do not perform well under such attacks, because nine agents (including six Byzantine agents and three regular agents) essentially use the training images of one digit, such that the models trained from this particular digit dominate. Therefore, the majority voting rules of ByRDiE-S and BRIDGE-S emphasize more on this particular digit, while ignore other digits relatively.
	
	\black{
		\noindent \textbf{Impact of Fraction of Byzantine Agents with Non-i.i.d. Data}. We further study the robustness of our proposed method with different fractions of Byzantine agents with Non-i.i.d. data. The data distributions and Byzantine attacks are the same as those in the above paragraph. In our proposed method, the penalty parameter is $\lambda=0.02$ and the step size is $\alpha^k = \frac{0.4}{\sqrt{k+1}}$. Not surprisingly, as shown in Figure \ref{Imob-noniid} the performance of all the methods degrades when the fraction of Byzantine agents increases. However, different to the i.i.d. case in Figure \ref{Imob-iid}, our proposed method shows much better robustness than ByRDiE-S and BRIDGE-S in the non-i.i.d. case.
	}
	
	\subsection{Numerical Experiments over Time-Varying Network}
	\label{sec:4b}
	
	Consider a static Erdos-Renyi networks consisting of $n=30$ agents, in which $b=3$ agents are Byzantine but the network of regular agents is connected. We generate two time-varying networks upon it. At every time $k$, every edge $e$ is randomly activated with probabilities $p_e = 0.01$ and $p_e = 0.005$, respectively. We consider the same-value attacks, where every Byzantine agent $j \in \mathcal{B}$ sends $z_j^k=c\mathbf{1}$ to its neighbors, with $\mathbf{1} \in \R^p$ being an all-one vector and $c=100$. The algorithm parameters are: (i) $\lambda=0.01$ and $\alpha^k=\frac{0.28}{\sqrt{k+1}}$ for the static network, (ii) $\lambda = 0.2$ and $\alpha^k = \frac{0.5}{\sqrt{k+1}}$ for the time-varying network with $p_e = 0.01$, and, (iii) $\lambda = 0.4$ and $\alpha^k = \frac{0.5}{\sqrt{k+1}}$ for the time-varying network with $p_e = 0.005$. Note that smaller $p_e$ means worse connectivity, which leads to smaller $\tilde{\sigma}_{\min}(\bar{A})$ as we have empirically observed from the experiments. Therefore, according to Theorem \ref{theorem4}, the critical value of $\lambda$ that guarantees consensus, given by $\lambda_0:=\frac{\sqrt{|\mathcal{R}|}}{\tilde{\sigma}_{\min}(\bar{A})} max_{i \in \mathcal{R}} \|\nabla\E[F(\tilde{x}^*,\xi_i)]+\nabla f_0(\tilde{x}^*)\|_\infty$, should be larger for smaller $p_e$.
	
\begin{figure}[H]
		\centering
		\includegraphics[scale=0.4]{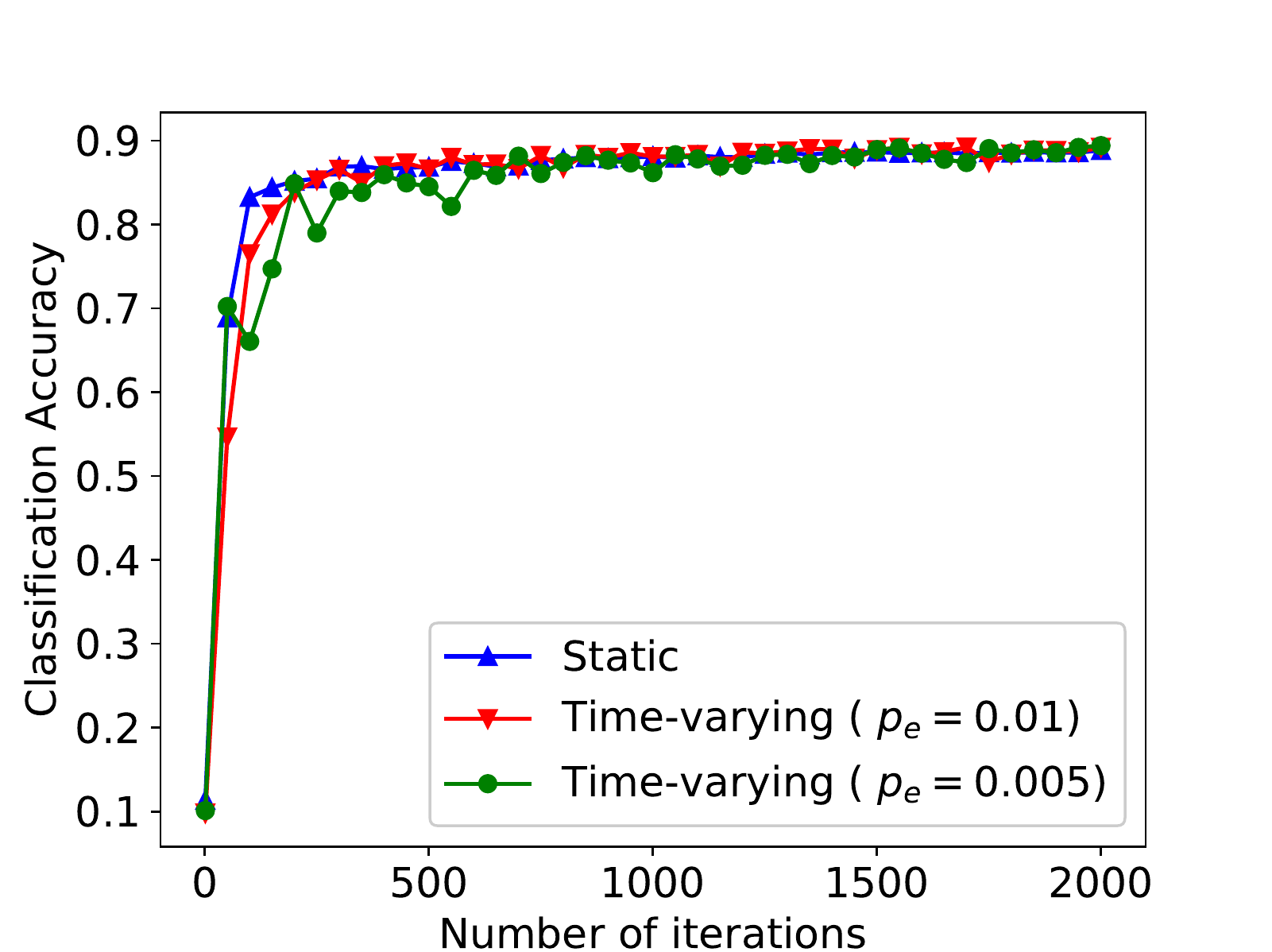} \hspace{-1.5em}
		\includegraphics[scale=0.4]{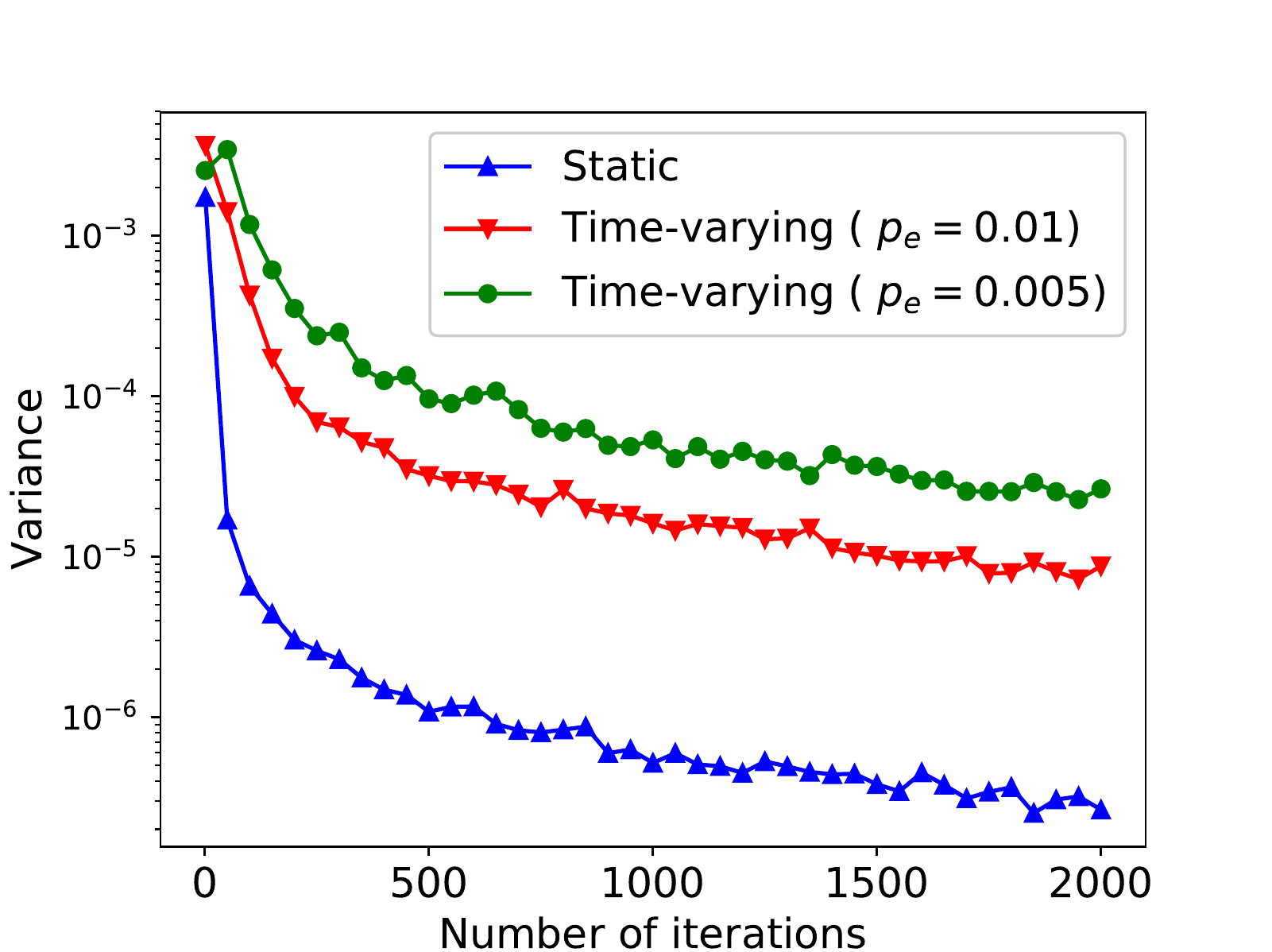}
		\caption{Classification accuracy and variance of regular agents' local models under same-value attacks over static and time-varying networks.} \label{COMP-SV}
	\end{figure}

	As shown in Figure \ref{COMP-SV}, our proposed method shows remarkable robustness to the Byzantine attacks even when the network is time-varying. The classification accuracies are almost the same over the static and the time-varying networks. The level of consensus degrades when $p_e$ decreases. This makes sense because in a less connected network, information diffusion is slower such that reaching consensus becomes more difficult.
	
	\vspace{-0.5em}
	
	\vspace{-0.5em}

	\section{Conclusion}
	\label{sec:5}
	
	In this paper, we deal with the Byzantine-robust decentralized stochastic optimization problem over static and time-varying networks. We introduce the TV norm penalty to handle Byzantine attacks and propose a stochastic subgradient method to solve the penalized problem. Theoretical analysis and numerical experiments demonstrate the robustness of our proposed method to Byzantine attacks, no matter the network topology is static or time-varying. In our future work, we will investigate the application of variance reduction techniques, which have been shown effective in Byzantine-robust centralized stochastic optimization \cite{Wu2019VR}, to the decentralized case. We will also consider TV norm penalties based on $\ell_2$ and $\ell_\infty$ norms \cite{RSA}.
	
	\vspace{1em}
	\noindent \textbf{Acknowledgement.} Qing Ling is supported in part by NSF China Grant 61973324, and Fundamental Research Funds for the Central Universities. A preliminary version of this paper has appeared in IEEE International Conference on Acoustics, Speech, and Signal Processing, Barcelona, Spain, May 4-8, 2020 \cite{Jie2020-icassp}.
	\vspace{1em}
	
	\appendix
	\section{Proof of Theorem \ref{theorem1}}
	\label{proof-theorem1}
	
	\noindent\textbf{Proof.} The optimal solution $x^*:=[x_i^*]$ of \eqref{eq6} satisfies the optimality condition
	\begin{align}\label{eq6-opt}
		\nabla \E[F(x_i^*, \xi_i)] + \lambda \sum_{j \in \mathcal{R}_i, i<j} s_{ij} - \lambda \sum_{j \in \mathcal{R}_i, i>j} \blue{s_{ji}} 
		+ \nabla f_0(x^*_i) = 0, \quad \forall i \in \mathcal{R},
	\end{align}
	where $s_{ij} = sign(x^*_i -x^*_j)\in[-1,1]^p$ for $i<j$. Note that $x^*$ is unique due to the strong convexity given by Assumption \ref{assumption1}. We will prove that $[\tilde{x}^*]$ satisfies \eqref{eq6-opt}, such that
	\begin{align}\label{eq6-opt-0}
		\nabla \E[F(\tilde{x}^*, \xi_i)] + \lambda \sum_{j \in \mathcal{R}_i, i<j} s_{ij} - \lambda \sum_{j \in \mathcal{R}_i, i>j} s_{ji}
		+ \nabla f_0(\tilde{x}^*) = 0, \quad \forall i \in \mathcal{R}.
	\end{align}

	For simplicity, define $v_i:=\nabla \E[F(\tilde{x}^*, \xi_i)]+\nabla f_0(\tilde{x}^*)$. Since \eqref{eq6-opt-0} can be decomposed element-wise, from now on we assume the variable dimension $p=1$ such that both $s_{ij}$ and $v_i$ are scalars. Then, \eqref{eq6-opt-0} can be rewritten as a system of linear equations
	\begin{align}\label{eq6-opt-1}
		\lambda A s + v = 0,
	\end{align}
	where $s =[s_{ij}]$ collects all the scalars $s_{ij}$ in order and $v=[v_i]$ collects all the scalars $v_i$ in order. Now the problem is equivalent to finding a vector $s$ whose elements are within $[-1,1]$, namely, $\|s\|_\infty \leq 1$, to satisfy \eqref{eq6-opt-1}.
	
	We first show that \eqref{eq6-opt-1} has at least one solution. To see so, observe that the rank of $A$ is $|\mathcal{R}|-1$ and the null space of the columns is spanned by the all-one vector $1_{|\mathcal{R}|} \in \mathbb{R}^{|\mathcal{R}|}$, because $(\mathcal{R},\mathcal{E}_{R})$ is bidirectionally connected according to Assumption \ref{assumption-static}. Meanwhile, according to the optimality condition of \eqref{eq4}, $\sum_{i\in \mathcal{R}} v_i = \sum_{i\in \mathcal{R}} \left( \nabla \E[F(\tilde{x}^*, \xi_i)]+\nabla f_0(\tilde{x}^*) \right)=0$. Therefore, the columns of $\lambda A$ and those of $[\lambda A, v]$ share the same null space and have the same rank. Consequently, we can find at least one solution to \eqref{eq6-opt-1}.
	
	We next find one solution to \eqref{eq6-opt-1} that satisfies $\|s\|_\infty \leq 1$. According to the above derivation, we also know that we can find at least one solution to $s'$ to $A s' + v = 0$. Among all the solutions to $A s' + v = 0$, we consider the least-squares solution given by $s'=-A^\dag v$, where $\dag$ denotes the pseudo inverse. This solution is bounded by $$\|s'\|_2 = \|A^\dag v\|_2 \leq \sigma_{\max}(A^\dag) \|v\|_2 \leq \frac{1}{\tilde{\sigma}_{\min}(A)} \|v\|_2,$$ where $\sigma_{\max}(\cdot)$ and $\tilde{\sigma}_{\min}(\cdot)$ denotes the largest and the smallest nonzero singular values, respectively. Since $\|s'\|_\infty \leq \|s'\|_2$ and $\|v\|_2 \leq \sqrt{|\mathcal{R}|} \|v\|_\infty$, we further have $$\|s'\|_\infty \leq \frac{\sqrt{|\mathcal{R}|}}{\tilde{\sigma}_{\min}(A)} \|v\|_\infty = \frac{\sqrt{|\mathcal{R}|}}{\tilde{\sigma}_{\min}(A)} max_{i \in \mathcal{R}} |v_i|.$$
	
	Then, we construct $s=\frac{1}{\lambda}s'$, which is a solution to \eqref{eq6-opt-1} because $A s' + v = 0$. In addition, we have $$\|s\|_\infty = \frac{1}{\lambda} \|s'\|_\infty \leq \frac{\sqrt{|\mathcal{R}|}}{\lambda \tilde{\sigma}_{\min}(A)} max_{i \in \mathcal{R}} |v_i|,$$ which is less than $1$ as long as $\lambda \geq \frac{\sqrt{|\mathcal{R}|}}{\tilde{\sigma}_{\min}(A)} max_{i \in \mathcal{R}} |v_i|$.
	
	Now we consider the variable dimension $p \geq 1$. For all the dimensions we construct vectors $s$ in this way such that $\|s\|_\infty \leq 1$ as long as $\lambda \geq \lambda_0 := \frac{\sqrt{|\mathcal{R}|}}{\tilde{\sigma}_{\min}(A)} max_{i \in \mathcal{R}} \|\nabla\E[F(\tilde{x}^*,\xi_i)]+\nabla f_0(\tilde{x}^*)\|_\infty$. This completes the proof.
	
	\section{Proof of Theorem \ref{theorem2}}
	\label{proof-theorem2}

	\noindent\textbf{Proof.}
	\textbf{Step 1.} \black{
		We first take conditional expectations given the variables up to time $k$, namely $\{x_i^l: l\le k, i \in \mathcal{R}\}$, and then take the expectation over $\{x_i^l: l\le k, i \in \mathcal{R}\}$. We simply denote $\E[\cdot | x_i^l: l\le k, i \in \mathcal{R}]$ as $\E_k[\cdot]$.} From the update \eqref{eq:subgradient} at every regular agent $i$, we have:
	\begin{align}\label{b1}
		&\black{\E_k}\|x_i^{k+1}-x_i^*\|^2\\
		=&\black{\E_k}\|x_i^{k}-x_i^*-\alpha^{k}\big(
		\nabla F(x_i^k,\xi_i^k)+ \lambda\sum_{j \in \mathcal{R}_i} sign(x_i^k-x_j^k)+\lambda\sum_{j \in \mathcal{B}_i}sign(x_i^k-z_j^k)+\nabla f_0(x_i^k)\big)\|^2 \nonumber\\
		=&\|x_i^k-x_i^*\|^2+(\alpha^{k})^2\black{\E_k}\|\nabla F(x_i^k,\xi_i^k)+ \lambda\sum_{j \in \mathcal{R}_i} sign(x_i^k-x_j^k)+\lambda\sum_{j \in \mathcal{B}_i}sign(x_i^k-z_j^k)+\nabla f_0(x_i^k)\|^2 \nonumber\\
		&-2\alpha^{k}\langle \black{\E_k}[\nabla F(x_i^k,\xi_i^k)]\ +\lambda\sum_{j \in \mathcal{R}_i} sign(x_i^k-x_j^k)+\nabla f_0(x_i^k), x_i^k-x_i^*\rangle
		-2\alpha^{k}\langle \lambda\sum_{j \in \mathcal{B}_i} sign(x_i^k-z_j^k), x_i^k-x_i^* \rangle. \nonumber
	\end{align}
	Below, we handle the terms at the right-hand side of \eqref{b1} one by one.
	
	For the second term at the right-hand side of \eqref{b1}, we have
	\begin{align} \label{b3}
		&\black{\E_k}\|\nabla F(x_i^k,\xi_i^k)+ \lambda\sum_{j \in \mathcal{R}_i} sign(x_i^k-x_j^k)+\lambda\sum_{j \in \mathcal{B}_i}sign(x_i^k-z_j^k)+\nabla f_0(x_i^k)\|^2 \\
		=&\black{\E_k}\|\nabla \black{\E_k}[F(x_i^k,\xi_i^k)]+ \lambda\sum_{j \in \mathcal{R}_i} sign(x_i^k-x_j^k)+\lambda\sum_{j \in \mathcal{B}_i}sign(x_i^k-z_j^k)+\nabla f_0(x_i^k) + \nabla F(x_i^k,\xi_i^k) - \nabla\black{\E_k}[F(x_i^k,\xi_i^k)]\|^2 \nonumber\\
		\leq&2\|\nabla \black{\E_k}[F(x_i^k,\xi_i^k)]+ \lambda\sum_{j \in \mathcal{R}_i} sign(x_i^k-x_j^k)+\lambda\sum_{j \in \mathcal{B}_i}sign(x_i^k-z_j^k)+\nabla f_0(x_i^k)\|^2 + 2\black{\E_k}\|F(x_i^k,\xi_i^k) - \nabla\black{\E_k}[F(x_i^k,\xi_i^k)]\|^2 \nonumber\\
		\leq&4\|\nabla \black{\E_k}[F(x_i^k,\xi_i^k)]+\lambda\sum_{j \in \mathcal{R}_i} sign(x_i^k-x_j^k)+\nabla f_0(x_i^k)\|^2+4\lambda^2\|\sum_{j \in \mathcal{B}_i}sign(x_i^k-z_j^k)\|^2 + 2\black{\E_k}\|F(x_i^k,\xi_i^k) - \nabla\black{\E_k}[F(x_i^k,\xi_i^k)]\|^2 \nonumber\\
		\leq&4\|\nabla \black{\E_k}[F(x_i^k,\xi_i^k)]+\lambda\sum_{j \in \mathcal{R}_i} sign(x_i^k-x_j^k)+\nabla f_0(x_i^k)\|^2+4\lambda^2|\mathcal{B}_i|^2p +2\delta_i^2, \nonumber
	\end{align}
	where the last inequality holds true because each element of the $p$-dimensional vector $sign(x_i^k-z_j^k)$ is within $[-1,1]$, and the variance is bounded by $\black{\E_k}\|F(x_i^k,\xi_i^k) - \nabla\black{\E_k}[F(x_i^k,\xi_i^k)]\|^2 \leq \delta_i^2$ stated in Assumption \ref{assumption3} and the independence of random $\xi_i^k$'s.
	\black{Recall the optimality condition of \eqref{eq6} shown in \eqref{eq6-opt}. For simplicity, we define $s_{ji} = -s_{ij}$ for $j>i$, then
		\begin{equation}\label{eq:opt}
			\nabla \E[F(x_i^*,\xi_i)] + \lambda\sum_{j \in \mathcal{R}_i} s_{ij} + \nabla f_0(x_i^*) = 0,
		\end{equation}
		which holds true for some $s_{ij} = sign(x_i^* - x_j^*)\in[-1,1]^p$. Plugging this equation} 
	into the first term at the right-hand side of \eqref{b3} and replacing $\xi_i$ by $\xi_i^k$, we have
	\black{
		\begin{align}\label{b4-tmp}
			&4\|\nabla \E_k[F(x_i^k,\xi_i^k)]+\lambda\sum_{j \in \mathcal{R}_i} sign(x_i^k-x_j^k)+\nabla f_0(x_i^k)\|^2 \\
			=&4\|\nabla \E_k[F(x_i^k,\xi_i^k)]+\lambda\sum_{j \in \mathcal{R}_i} sign(x_i^k-x_j^k)+\nabla f_0(x_i^k)-\nabla \E_k[F(x_i^*,\xi_i^k)]-\lambda\sum_{j \in \mathcal{R}_i} \black{s_{ij}}-\nabla f_0(x_i^*) \|^2 \nonumber\\
			\leq& 8\|\nabla \E_k[F(x_i^k,\xi_i^k)] + \nabla f_0(x_i^k) - \nabla \E_k[F(x_i^*,\xi_i^k)] - \nabla f_0(x_i^*)\|^2 + 8\lambda^2\|\sum_{j \in \mathcal{R}_i}sign(x_i^k-x_j^k)-\sum_{j \in \mathcal{R}_i} s_{ij}  \|^2 \nonumber \\
			\leq& 8\|\nabla \E_k[F(x_i^k,\xi_i^k)] + \nabla f_0(x_i^k) - \nabla \E_k[F(x_i^*,\xi_i^k)] - \nabla f_0(x_i^*)\|^2 + 32\lambda^2|\mathcal{R}_i|^2p \nonumber,
		\end{align}
	}
	where the first inequality is due to $(a+b)^2\leq2a^2+2b^2$. Combining \eqref{b3} and \eqref{b4-tmp}, we have
	\begin{align}\label{b4}
		&\black{\E_k}\|\nabla F(x_i^k,\xi_i^k)+ \lambda\sum_{j \in \mathcal{R}_i} sign(x_i^k-x_j^k)+\lambda\sum_{j \in \mathcal{B}_i}sign(x_i^k-z_j^k)+\nabla f_0(x_i^k)\|^2 \\
		\leq& \black{8\|\nabla \black{\E_k}[F(x_i^k,\xi_i^k)] + \nabla f_0(x_i^k) - \nabla \black{\E_k}[F(x_i^*,\xi_i^k)] - \nabla f_0(x_i^*)\|^2 + 32}
		\lambda^2|\mathcal{R}_i|^2p+4\lambda^2|\mathcal{B}_i|^2p +2\delta_i^2. \nonumber
	\end{align}

	For the third term at the right-hand side of \eqref{b1}, noticing that $x_i^k$ is independent with $\xi_i^k$ such that
	\begin{align}\label{b2}
		\langle \black{\E_k}[\nabla F(x_i^k,\xi_i^k)], x_i^k-x_i^* \rangle = \langle \nabla \black{\E_{\xi_i^k}}[F(x_i^k,\xi_i^k)], x_i^k-x_i^* \rangle,
	\end{align}
	we have
	\begin{align}\label{b2-tmp}
		& -2\langle \black{\E_k}[\nabla F(x_i^k,\xi_i^k)]+\lambda\sum_{j \in \mathcal{R}_i} sign(x_i^k-x_j^k)+\nabla f_0(x_i^k), x_i^k-x_i^*\rangle \\
		= & -2 \langle \nabla \black{\E_{\xi_i^k}}[F(x_i^k,\xi_i^k)]+\lambda\sum_{j \in \mathcal{R}_i} sign(x_i^k-x_j^k)+\nabla f_0(x_i^k), x_i^k-x_i^*\rangle. \nonumber
	\end{align}
	Plugging the optimality condition \eqref{eq:opt} into \eqref{b2-tmp} , we have
	\begin{align}\label{b5}
		&-2\langle \nabla \black{\E_{\xi_i^k}}[F(x_i^k,\xi_i^k)]+\lambda\sum_{j \in \mathcal{R}_i} sign(x_i^k-x_j^k)+\nabla f_0(x_i^k), x_i^k-x_i^*\rangle\\
		=&-2\langle\nabla \black{\E_{\xi_i^k}}[F(x_i^k,\xi_i^k)] + \nabla f_0(x_i^k) - \nabla \black{\E_{\xi_i^k}}[F(x_i^*,\xi_i^k)] - \nabla f_0(x_i^*), x_i^k-x_i^*\rangle -2\langle\lambda\sum_{j \in \mathcal{R}_i}sign(x_i^k-x_j^k)-\lambda\sum_{j \in \mathcal{R}_i}\black{s_{ij}} , x_i^k-x_i^*\rangle \nonumber
	\end{align}
	Since we assume that the functions $\black{\E_{\xi_i^k}}[F(x_i,\black{\xi_i^k})] + f_0(x_i)$ are strongly convex and have Lipschitz continuous gradients (cf. Assumption \ref{assumption1} and \ref{assumption2}), by \cite{Nesterov2004} we have
	\black{
		\begin{align}\label{b6}
			& -2\langle\nabla \black{\E_{\xi_i^k}}[F(x_i^k,\xi_i^k)] + \nabla f_0(x_i^k) - \nabla \black{\E_{\xi_i^k}}[F(x_i^*,\xi_i^k)] - \nabla f_0(x_i^*), x_i^k-x_i^*\rangle \\
			\leq & -\frac{2u_iL_i}{u_i+L_i}  \|x_i^k-x_i^*\|^2 -\frac{2}{u_i+L_i}\|\nabla \black{\E_{\xi_i^k}}[F(x_i^k,\xi_i^k)] + \nabla f_0(x_i^K) - \nabla \black{\E_{\xi_i^k}}[F(x_i^*,\xi_i^k)]^2 - \nabla f_0(x_i^*) \|^2. \nonumber
		\end{align}
	}
	Substituting \eqref{b5}, \eqref{b6} into \eqref{b2-tmp}, we have
	\black{
		\begin{align}\label{b8}
			&-2\langle \E_k[\nabla F(x_i^k,\xi_i^k)]+\lambda\sum_{j \in \mathcal{R}_i} sign(x_i^k-x_j^k)+\nabla f_0(x_i^k), x_i^k-x_i^*\rangle\\
			\leq&-2\langle\lambda\sum_{j \in \mathcal{R}_i}sign(x_i^k-x_j^k)-\lambda\sum_{j \in \mathcal{R}_i}\black{s_{ij}} , x_i^k-x_i^*\rangle \nonumber
			-(\frac{2u_iL_i}{u_i+L_i}) \  \|x_i^k-x_i^*\|^2 \nonumber \\
			&-\frac{2}{u_i+L_i} \ \|\nabla \E_{\xi_i^k}[F(x_i^k,\xi_i^k)] + \nabla f_0(x_i^k) - \nabla \E_{\xi_i^k}[F(x_i^*,\xi_i^k)] - \nabla f_0(x_i^*)\|^2 . \nonumber
		\end{align}
	}
	
	For the last term at the right-hand side of \eqref{b1}, it holds for any $\epsilon>0$ that
	\black{
		\begin{align}\label{b9}
			&-2\langle \lambda\sum_{j \in \mathcal{B}_i} sign(x_i^k-z_j^k), x_i^k-x_i^*   \rangle \\
			\leq &\epsilon \|x_i^k-x_i^*\|^2 + \frac{\lambda^2}{\epsilon} \|\sum_{j \in \mathcal{B}_i} sign(x_i^k-z_j^k)\|^2 \nonumber\\
			\leq& \epsilon \|x_i^k-x_i^*\|^2  + \frac{\lambda^2|\mathcal{B}_i|^2p}{\epsilon}. \nonumber
		\end{align}
	}
	
	Substituting \eqref{b4}, \eqref{b8} and \eqref{b9} into \eqref{b1} and combining the terms, we have
	\black{\begin{align}\label{b10-tmp}
			&\black{\E_k}\|x_i^{k+1}-x_i^*\|^2\\
			\leq&\big(1-\alpha^{k}(\frac{2u_iL_i}{u_i+L_i}-\epsilon)\big)\|x_i^k-x_i^*\|^2 + (\alpha^{k})^2(32\lambda^2|\mathcal{R}_i|^2p+4\lambda^2|\mathcal{B}_i|^2p+2\delta_i^2)+\alpha^{k} \frac{\lambda^2|\mathcal{B}_i|^2p}{\epsilon} \nonumber \\
			&-2\alpha^{k}\langle\lambda\sum_{j \in \mathcal{R}_i}sign(x_i^k-x_j^k)-\lambda\sum_{j \in \mathcal{R}_i}\black{s_{ij}} , x_i^k-x_i^*\rangle \nonumber \\
			&-2\alpha^{k} (\frac{1}{u_i+L_i}-4\alpha^{k}) \|\nabla \black{\E_{\xi_i^k}}[F(x_i^k,\xi_i^k)] + \nabla f_0(x_i^k) - \nabla \black{\E_{\xi_i^k}}[F(x_i^*,\xi_i^k)] - \nabla f_0(x_i^*)\|^2. \nonumber
	\end{align}}
	Due to the step size rule, we have \black{$\frac{1}{u_i+L_i}-4\alpha^{k} \geq 0$} and hence can drop the last term of at the right-hand side of \eqref{b10-tmp}. Also noticing the definition of $\eta$, we rewrite \eqref{b10-tmp} into
	\begin{align}\label{b10}
		&\black{\E_k}\|x_i^{k+1}-x_i^*\|^2 \\
		\leq&\big(1-\eta\alpha^{k}\big)\|x_i^k-x_i^*\|^2 +(\alpha^{k})^2(\black{32}\lambda^2|\mathcal{R}_i|^2p+4\lambda^2|\mathcal{B}_i|^2p+2\delta_i^2)+\alpha^{k} \frac{\lambda^2|\mathcal{B}_i|^2p}{\epsilon} \nonumber \\
		&-2\alpha^{k}\langle\lambda\sum_{j \in \mathcal{R}_i}sign(x_i^k-x_j^k)-\lambda\sum_{j \in \mathcal{R}_i}\black{s_{ij}} , x_i^k-x_i^*\rangle. \nonumber
	\end{align}

	\noindent\textbf{Step 2.} Here we define $I_p(x)=\frac{\lambda}{2}\sum_{i \in \mathcal{R}}\sum_{j \in \mathcal{R}_i}\|x_i-x_j\|_1$. Since $I_p(x)$ is convex, we have
	\begin{align}\label{b.2.1}
		&\langle \partial_x I_p(x^k)- \partial_x I_p(x^*), x^k-x^* \rangle \\
		=&\sum_{i \in \mathcal{R}} \langle\lambda\sum_{j \in \mathcal{R}_i}sign(x_i^k-x_j^k)-\lambda\sum_{j \in \mathcal{R}_i}\black{s_{ij}} , x_i^k-x_i^*\rangle \geq 0. \nonumber
	\end{align}
	Summing up \eqref{b10} over all regular agents $i \in \mathcal{R}$ and adding to \eqref{b.2.1}, we have
	\begin{align} \label{b.2.2}
		&\black{\E_k}\|x^{k+1}-x^*\|^2\\
		\leq&\big(1-\eta\alpha^{k}\big)\|x^k-x^*\|^2 +(\alpha^{k})^2 \ \sum_{i \in \mathcal{R}} \ (\black{32}\lambda^2|\mathcal{R}_i|^2p+4\lambda^2|\mathcal{B}_i|^2p+2\delta_i^2) +\alpha^{k}\sum_{i \in \mathcal{R}}\frac{\lambda^2|\mathcal{B}_i|^2p}{\epsilon}. \nonumber
	\end{align}
	\black{Taking the full expectation of \eqref{b.2.2}, then we get
		\begin{align}
			\E \|x^{k+1 } - x^*\|^2 = \E[\E_k \|x^{k+1 } - x^*\|^2                                                                                                                                                                                                                                                                                                                                                                                                                                                                                                                                                                                                                                                                                                                                                                                                                                                                                                                                                                                                                                                                                                                                                                                                                                                                                                                                                                                                                                                                                                                                                                                                                                                                                                                                                                                                                                                                                                                                                                                                                                                                                                                                                                                                                                                                                                                                                                                                                                                                                                                                                                                                                                                                                                                                                                                                                                                                                                                                                                                                                                                                                                                                                                                                                                                                                                                                                                                                                                                                                                                                                                                                                                                                                                                                                                                                                                                                                                                                                                                                                                                                                                                                                                                                                                                                                                                                                                                                                                                                                                                                                                                                                                                                                                                                                                                                                                                                                                                                                                                                                                                                                                                                                                                                                                                                                                                                                                                                                                                                                                                                                                                                                                                                                                                                                                                                                                                                                                                                                                                                                                                                                                                                                                                                                                                                                                                                                                                                                                                                                                                                                                                                                                                                                                                                                                                                                                                                                                                                                                                                                                                                                                                                                                                                                                                                                                                                                                                                                                                                                                                                                                                                                                                                                                                                                                                                                                                                                                                                                                                                                                                                                                                                                                                                                                                                                                                                                                                                                                                                                                                                                                                                                                                                                                                                                                                                                                                                                                                                                                                                                                                                                                                                                                                                                                                                                                                                                                                                                                                                                                                                                                                                                                                                                                                                                                                                                                                                                                                                                                                                                                                                                                                                                                                                                                                                                                                                                                                                                                                                                                                                                                                                                                                                                                                                                                                                                                                                                                                                                                                                                                                                                                                                                                                                                                                                                                                                                                                                                                                                                                                                                                                                                                                                                                                                                                                                                                                                                                                                                                                                                                                                                                                                                                                                                                                                                                                                                                                                                                                                                                                                                                                                                                                                                                                                                                                                                                                                                                                                                                                                                                                                                                                                                                                                                                                                                                                                                                                                                                                                                                                                                                                                                                                                                                                                                                                                                                                                                                                                                                                                                                                                                                                                                                                                                                                                                                                                                                                                                                                                                                                                                                                                                                                                                                                                                                                                                                                                                                                                                                                                                                                                                                                                                                                                                                                                                                                                                                                                                                                                                                                                                                                                                                                                                                                                                                                                                                                                                                                                                                                                                                                                                                                                                                                                                                                                                                                                                                                                                                                                                                                                                                                                                                                                                                                                                                                                                                                                                                                                                                                                                                                                                                                                                                                                                                                                                                                                                                                                                                                                                                                                                                                                                                                                                                                                                                                                                                                                                                                                                                                                                                                                                                                                                                                                                                                                                                                                                                                                                                                                                                                                                                                                                                                                                                                                                                                                                                                                                                                                                                                                                                                                                                                                                                                                                                                                                                                                                                                                                                                                                                                                                                                                                                                                                                                                                                                                                                                                                                                                                                                                                                                                                                                                                                                                                                                                                                                                                                                                                                                                                                                                                                                                                                                                                                                                                                                                                                                                                                                                                                                                                                                                                                                                                                                                                                                                                                                                                                                                                                                                                                                                                                                                                                                                                                                                                                                                                                                                                                                                                                                                                                                                                                                                                                                                                                                                                                                                                                                                                                                                                                                                                                                                                                                                                                                                                                                                                                                                                                                                                                                                                                                                                                                                                                                                                                                                                                                                                                                                                                                                                                                                                                                                                                                                                                                                                                                                                                                                                                                                                                                                                                                                                                                                                                                                                                                                                                                                                                                                                                                                                                                                                                                                                                                                                                                                                                                                                                                                                                                                                                                                                                                                                                                                                                                                                                                                                                                                                                                                                                                                                                                                                                                                                                                                                                                                                                                                                                                                                                                                                                                                                                                                                                                                                                                                                                                                                                                                                                                                                                                                                                                                                                                                                                                                                                                                                                                                                                                                                                                                                                                                                                                                                                                                                                                                                                                                                                                                                                                                                                                                                                                                                                                                                                                                                                                                                                                                                                                                                                                                                                                                                                                                                                                                                                                                                                                                                                                                                                                                                                                                                                                                                                                                                                                                                                                                                                                                                                                                                                                                                                                                                                                                                                                                                                                                                                                                                                                                                                                                                                                                                                                                                                                                                                                                                                                                                                                                                                                                                                                                                                                                                                                                                                                                                                                                                                                                                                                                                                                                                                                                                                                                                                                                                                                                                                                                                                                                                                                                                                                                                                                                                                                                                                                                                                                                                                                                                                                                                                                                                                                                                                                                                                                                                                                                                                                                                                                                                                                                                                                                                                                                                                                                                                                                                                                                                                                                                                                                                                                                                                                                                                                                                                                                                                                                                                                                                                                                                                                                                                                                                                                                                                                                                                                                                                                                                                                                                                                                                                                                                                                                                                                                                                                                                                                                                                                                                                                                                                                                                                                                                                                                                                                                                                                                                                                                                                                                                                                                                                                                                                                                                                                                                                                                                                                                                                                                                                                                                                                                                                                                                                                                                                                                                                                                                                                                                                                                                                                                                                                                                                                                                                                                                                                                                                                                                                                                                                                                                                                                                                                                                                                                                                                                                                                                                                                                                                                                                                                                                                                                                                                                                                                                                                                                                                                                                                                                                                                                                                                                                                                                                                                                                                                                                                                                                                                                                                                                                                                                                                                                                                                                                                                                                                                                                                                                                         ] \leq (1 - \eta\alpha^k) \E \|x^k - x^*\|^2 + (\alpha^k)^2 \Delta_0 + \alpha^k \Delta_2,
		\end{align}
		}
	where the constants $\Delta_0$ and $\Delta_2$ are defined as
	\begin{align} \nonumber
		\Delta_0 = \sum\limits_{i \in \mathcal{R}} \big(&\black{32} \lambda^2|\mathcal{R}_i|^2p+4 \lambda^2|\mathcal{B}_i|^2p +2\delta_i^2\big), \quad \Delta_2 = \sum\limits_{i \in \mathcal{R}}\frac{\lambda^2|\mathcal{B}_i|^2p}{\epsilon}.
	\end{align}

	\noindent\textbf{Step 3.} According to the step size rule $\alpha^{k}=\min\{\underline{\alpha},\frac{\overline{\alpha}}{k+1}\}$, there exists a smallest integer $k_0$ satisfying $\underline{\alpha} \geq \frac{\overline{\alpha}}{k_0+1}$ such that $\alpha^{k}=\underline{\alpha}$ when $k<k_0$ and $\alpha^{k}=\frac{\overline{\alpha}}{k+1}$ when $k \geq k_0$. Then for all $k<k_0$, \eqref{b.2.2} becomes
	\begin{align} \label{b.3.1}
		&\E\|x^{k+1}-x^*\|^2 \leq  \big(1-\eta\underline{\alpha}\big)\E\|x^k-x^*\|^2 + (\underline{\alpha})^2\Delta_0 +\underline{\alpha}\Delta_2, \quad\forall k <k_0.
	\end{align}
	By the definitions of $\eta$ and $\underline{\alpha}$, $\eta\underline{\alpha}\in(0,1)$. Applying telescopic cancellation to \eqref{b.3.1} through time 0 to $k<k_0$ yields
	\begin{align} \label{b.3.2}
		&\E\|x^{k+1} - x^*\|^2 \leq (1-\eta\underline{\alpha})^k \E\|x^0 - x^*\|^2  + \frac{1}{\eta}(\underline{\alpha}\Delta_0+\Delta_2), \quad \forall{k}<k_0.
	\end{align}

	For all $k\geq k_0$, \eqref{b.2.2} becomes
	\begin{align} \label{b.3.3}
		&\E\|x^{k+1}-x^*\|^2 \leq \big(1-\frac{\eta\overline{\alpha}}{k+1}\big)\E\|x^k-x^*\|^2 + \frac{(\overline{\alpha})^2\Delta_0}{(k+1)^2} +\frac{\overline{\alpha}\Delta_2}{k+1}, \quad\forall k \geq k_0.
	\end{align}
	Note that $1-\frac{\eta\overline{\alpha}}{k+1} \in (0,1)$ when $k \geq k_0$. Below, we use induction to prove
	\begin{align} \label{b.3.4}
		\E\|x^{k+1} - x^*\|^2 \leq \frac{\Delta_1}{k+1}+\overline{\alpha}\Delta_2. \quad \forall{k} \geq k_0,
	\end{align}
	where
	\begin{align}
		\Delta_1=\max\{\frac{\overline{\alpha}^2\Delta_0}{\eta\overline{\alpha}-1},(k_0+1)\E\|x^{k_0}-x^*\|^2+\frac{\overline{\alpha}^2\Delta_0}{k_0+1}\}. \nonumber
	\end{align}
	When $k=k_0$, by \eqref{b.3.3}, we have
	\begin{align}\label{b.3.5}
		&\E\|x^{k_0+1}-x^*\|^2 \\
		\leq &\big(1-\frac{\eta\overline{\alpha}}{k_0+1}\big)\E\|x^{k_0}-x^*\|^2 + \frac{(\overline{\alpha})^2\Delta_0}{(k_0+1)^2} +\frac{\overline{\alpha}\Delta_2}{k_0+1}  \nonumber \\
		\leq&\E\|x^{k_0}-x^*\|^2 + \frac{(\overline{\alpha})^2\Delta_0}{(k_0+1)^2} +\frac{\overline{\alpha}\Delta_2}{k_0+1}  \nonumber \\
		\leq&\frac{\Delta_1}{k_0+1}+\overline{\alpha}\Delta_2 \nonumber.
	\end{align}
	Now suppose that \eqref{b.3.4} is true when $k=k'>k_0$, such that
	\begin{align}\label{b.3.6}
		\E\|x^{k'+1}-x^*\|^2 \leq \frac{\Delta_1}{k'+1} + \overline{\alpha}\Delta_2.
	\end{align}
	When $k=k'+1$, because $k'+1>k_0$ according to \eqref{b.3.3} we have
	\begin{align}\label{b.3.7}
		&\E\|x^{k'+2}-x^*\|^2 \leq \big(1-\frac{\eta\overline{\alpha}}{k'+2}\big)\E\|x^{k'+1}-x^*\|^2 + \frac{(\overline{\alpha})^2\Delta_0}{(k'+2)^2} +\frac{\overline{\alpha}\Delta_2}{k'+2}.
	\end{align}
	Substituting \eqref{b.3.6} into \eqref{b.3.7}, we have:
	\begin{align}
		&\E\|x^{k'+2}-x^*\|^2 \\
		\leq&\big(1-\frac{\eta\overline{\alpha}}{k'+2}\big) \frac{\Delta_1}{k'+1} + \frac{(\overline{\alpha})^2\Delta_0}{(k'+2)^2} + \overline{\alpha}\Delta_2+ (1-\eta\overline{\alpha})\frac{\overline{\alpha}\Delta_2}{k'+2} \nonumber\\
		\overset{(a)}{\leq}&\big(1-\frac{\eta\overline{\alpha}}{k'+2}\big) \frac{\Delta_1}{k'+1} + \frac{(\overline{\alpha})^2\Delta_0}{(k'+2)^2} + \overline{\alpha}\Delta_2 \nonumber \\
		\overset{(b)}{\leq}&\big(1-\frac{\eta\overline{\alpha}}{k'+2}\big) \frac{\Delta_1}{k'+1} + \frac{(\eta\overline{\alpha}-1)\Delta_1}{(k'+2)^2} + \overline{\alpha}\Delta_2 \nonumber \\
		\leq&\big(1-\frac{\eta\overline{\alpha}}{k'+2}\big) \frac{\Delta_1}{k'+1} + \frac{(\eta\overline{\alpha}-1)\Delta_1}{(k'+1)(k'+2)} + \overline{\alpha}\Delta_2 \nonumber \\ \leq&\frac{\Delta_1}{k'+2}+\overline{\alpha}\Delta_2. \nonumber
	\end{align}
	where $(a)$ uses the fact that $\overline{\alpha}>\frac{1}{\eta}$, and $(b)$ follows from $\Delta_1 \geq \frac{\overline{\alpha}^2\Delta_0}{\eta\overline{\alpha}-1}$. This completes the induction and the entire proof.
	
	\section{Proof of Theorem \ref{theorem3}}
	\label{proof-theorem3}
	
	\noindent\textbf{Proof.} When $\lambda\geq\lambda_0$, combining Theorems \ref{theorem1} and \ref{theorem2}, we directly reach \eqref{th3-1}. When $0<\lambda<\lambda_0$, we have
	\begin{align}
		\E[\|x^{k+1}-[\tilde{x}^*]\|^2]\leq 2\E[\|x^{k+1}-x^*\|^2] +2 \E[\|x^*-[\tilde{x}^*]\|^2].
	\end{align}
	Because $\E\|x^{k+1}-x^*\|^2 \leq \frac{\Delta_1}{k+1}+\overline{\alpha}\Delta_2 \ $ and $\ \E[\|x^*-[\tilde{x}^*]\|^2]\leq\Delta_3$, we reach \eqref{th3-2} and complete the proof.
	
	\section{Proof of Theorem \ref{theorem4}}
	\label{proof-theorem4}
	
	\noindent\textbf{Proof.} The optimal solution $x^*:=[x_i^*]$ of \eqref{tvn-eq2} satisfies the optimality condition that for any $i\in\mathcal{R}$,
	\black{
		\begin{align}\label{tvn-eq2-opt}
			0 &= \nabla \E_{\xi_i}[F(x_i^*, \xi_i)] +  \E_{\bar\zeta} \bigg[\lambda \sum_{j \in \mathcal{R}_i(\bar\zeta), i<j} \bar{s}_{ij} - \lambda \sum_{j \in \mathcal{R}_i(\bar\zeta), i>j} \bar{s}_{ji}\bigg] + \nabla f_0(x^*_i),\\
			& = \nabla \E_{\xi_i}[F(x_i^*, \xi_i)] +  \lambda \sum_{e = (i,j) \in \bar{\mathcal{E}}_R, i<j} \bar{a}_{e} \bar{s}_{ij} - \lambda \sum_{e = (i,j) \in \bar{\mathcal{E}}_R, i>j} \bar{a}_{e} \bar{s}_{ji} + \nabla f_0(x^*_i), \label{tvn-eq2-opt-tmp}
		\end{align}
		for some $\bar{s}_{ij} = sign(x^*_i -x^*_j)\in[-1,1]^p$.} 
	Note that $x^*$ is unique due to the strong convexity given by Assumption \ref{assumption1}. 
	%
	%
	We will prove that $[\tilde{x}^*]$ satisfies \eqref{tvn-eq2-opt-tmp}, such that
	\begin{align}\label{tvn-eq2-opt-0}
		\nabla \black{\E_{\xi_i}}[F(\tilde{x}^*, \xi_i)] +  \lambda \sum_{e = (i,j) \in \bar{\mathcal{E}}_R, i<j} \black{\bar{a}_{e} \bar{s}_{ij}} - \lambda \sum_{e = (i,j) \in \bar{\mathcal{E}}_R, i>j} \black{\bar{a}_{e} \bar{s}_{ji}} + \nabla f_0(\tilde{x}^*) = 0, \quad \forall i \in \mathcal{R}.
	\end{align}

	Define $v_i:=\nabla \E[F(\tilde{x}^*, \xi_i)]+\nabla f_0(\tilde{x}^*)$. Since \eqref{tvn-eq2-opt-0} can be decomposed element-wise, we start from assuming the variable dimension $p=1$ such that both \black{$\bar{s}_{ij}$} and $v_i$ are scalars, and then extend to the high-dimensional case. By the definition of $\bar{A}$, \eqref{tvn-eq2-opt-0} can be rewritten as a system of linear equations
	\begin{align}\label{tvn-eq2-opt-1}
		\lambda \bar{A} \black{\bar{s}} + v = 0,
	\end{align}
	where $\black{\bar{s}} =[\black{\bar{s}_{ij}}]$ collects all the scalars $\black{\bar{s}_{ij}}$ in order, $v=[v_i]$ collects all the scalars $v_i$ in order. Now the problem is equivalent to finding a vector $\black{\bar{s}}$ whose elements are within $[-1,1]$, namely, $\|\black{\bar{s}}\|_\infty \leq 1$, to satisfy \eqref{tvn-eq2-opt-1}. The rest of the proof is the same as that of Theorem \ref{theorem1}, only replacing $A$ by $\bar{A}$.
	
	\section{Proof of Theorem \ref{theorem5}}
	\label{proof-theorem5}
	
	\noindent\textbf{Proof}. \black{In the proof, we focus on the sequence $\{\zeta^k\}$ satisfying $\lim_{K\rightarrow\infty}
		\frac1{K+1}\sum_{k=0}^K\zeta_e^k = \bar{a}_e$ for all $e$, which happens with probability $1$ according to Assumption \ref{assumption-timevary}. Same as the proof of Theorem \ref{theorem2}, we use the simplified notation $\E_k$ to take the conditional expectation given $\{x_i^l: l\le k, i \in \mathcal{R}\}$, and impose no restrictions on the attacks.}
	
	\noindent\textbf{Step 1.} From the update \eqref{tvn-eq6} at every regular agent $i$, we have
	\begin{align}\label{e1}
		&\black{\E_k}\|x_i^{k+1}-x_i^*\|^2\\
		=&\black{\E_k}\|x_i^{k}-x_i^*-\alpha^{k}\big(
		\nabla F(x_i^k,\xi_i^k)+ \lambda\sum_{j \in \mathcal{R}_i^k} sign(x_i^k-x_j^k)+\lambda\sum_{j \in \mathcal{B}_i^k}sign(x_i^k-z_j^k)+\nabla f_0(x_i^k)\big)\|^2 \nonumber\\
		=&\|x_i^k-x_i^*\|^2+(\alpha^{k})^2\black{\E_k}\|\nabla F(x_i^k,\xi_i^k)+ \lambda\sum_{j \in \mathcal{R}_i^k} sign(x_i^k-x_j^k)+\lambda\sum_{j \in \mathcal{B}_i^k}sign(x_i^k-z_j^k)+\nabla f_0(x_i^k)\|^2 \nonumber\\
		&-2\alpha^{k}\black{\E_k}\langle \nabla F(x_i^k,\xi_i^k)\ +\lambda\sum_{j \in \mathcal{R}_i^k} sign(x_i^k-x_j^k) +\nabla f_0(x_i^k), x_i^k-x	_i^*\rangle -2\alpha^{k}\langle \lambda\sum_{j \in \mathcal{B}_i^k} sign(x_i^k-z_j^k), x_i^k-x_i^* \rangle. \nonumber
	\end{align}
	Below, we handle the terms at the right-hand side of \eqref{e1} one by one.
	
	For the second term at the right-hand side of \eqref{e1}, we have
	\begin{align} \label{e3}
		&\black{\E_k}\|\nabla F(x_i^k,\xi_i^k)+ \lambda\sum_{j \in \mathcal{R}_i^k} sign(x_i^k-x_j^k)+\lambda\sum_{j \in \mathcal{B}_i^k}sign(x_i^k-z_j^k)+\nabla f_0(x_i^k)\|^2 \\
		=&\black{\E_k}\|\nabla \black{\E_{\xi_i^k}}[F(x_i^k,\xi_i^k)]+ \lambda\sum_{j \in \mathcal{R}_i^k} sign(x_i^k-x_j^k)+\lambda\sum_{j \in \mathcal{B}_i^k}sign(x_i^k-z_j^k)+\nabla f_0(x_i^k) + \nabla F(x_i^k,\xi_i^k) - \nabla\black{\E_{\xi_i^k}}[F(x_i^k,\xi_i^k)]\|^2 \nonumber\\
		\leq&4\|\nabla \black{\E_{\xi_i^k}}[F(x_i^k,\xi_i^k)]+\lambda\sum_{j \in \mathcal{R}_i^k} sign(x_i^k-x_j^k)+\nabla f_0(x_i^k)\|^2+4\lambda^2\|\sum_{j \in \mathcal{B}_i^k}sign(x_i^k-z_j^k)\|^2 + 2\black{\E_{\xi_i^k}}\|F(x_i^k,\xi_i^k) - \nabla\black{\E_{\xi_i^k}}[F(x_i^k,\xi_i^k)]\|^2 \nonumber \\
		\leq& 4\|\nabla \black{\E_{\xi_i^k}}[F(x_i^k,\xi_i^k)]+\lambda\sum_{j \in \mathcal{R}_i^k} sign(x_i^k-x_j^k)+\nabla f_0(x_i^k)\|^2+4\lambda^2|\mathcal{B}_i^k|^2p +2\delta_i^2 \nonumber\\
		\leq& 4\|\nabla \black{\E_{\xi_i^k}}[F(x_i^k,\xi_i^k)]+\lambda\sum_{j \in \mathcal{R}_i^k} sign(x_i^k-x_j^k)+\nabla f_0(x_i^k)\|^2+4\lambda^2|\mathcal{B}_i|^2p +2\delta_i^2, \nonumber
	\end{align}
	where the second equality holds true because that each element of the $p$-dimensional vector $sign(x_i^k-z_j^k)$ is within $[-1,1]$, and that the variance is bounded by $\E\|F(x_i^k,\xi_i^k) - \nabla\E[F(x_i^k,\xi_i^k)]\|^2 \leq \delta_i^2$ stated in Assumption \ref{assumption3}.
	
	\black{Let $\bar{s}_{ij}=-\bar{s}_{ji}$ for $i>j$. Then the optimality condition \eqref{tvn-eq2-opt} becomes $0 = \nabla \E_{\xi_i}[F(x_i^*,\xi_i)] + \lambda \ \E_{\bar\zeta}[\sum_{j \in \mathcal{R}_i(\bar\zeta)} \bar{s}_{ij}]+\nabla f_0(x_i^*)$. Plugging it} into the first term at the right-hand side of \eqref{e3} and replacing $\xi_i$ by $\xi_i^k$, we have
	\begin{align}\label{e4}
		&4\|\nabla\black{\E_{\xi_i^k}}F(x_i^k,\xi_i^k) + \lambda\sum_{j \in \mathcal{R}_i^k} sign(x_i^k-x_j^k) + \nabla f_0(x_i^k)\|^2 \\
		=&4\|\nabla\black{\E_{\xi_i^k}}[F(x_i^k,\xi_i^k)]+\lambda\sum_{j \in \mathcal{R}_i^k} sign(x_i^k-x_j^k)+\nabla f_0(x_i^k)-\nabla \black{\E_{\xi_i^k}}[F(x_i^*,\xi_i^k)]-\lambda \black{\E_{\bar\zeta}}[\sum_{j \in \mathcal{R}_i(\black{\bar\zeta})} \black{\bar{s}_{ij}}]-\nabla f_0(x_i^*) \|^2 \nonumber\\
		\leq& \black{8\|\nabla\E_{\xi_i^k}[F(x_i^k,\xi_i^k)] +\nabla f_0(x_i^k) - \nabla \E_{\xi_i^k}[F(x_i^*,\xi_i^k)] -\nabla f_0(x_i^*)\|^2
			+ 8\lambda^2\|\sum_{j \in \mathcal{R}_i^k}sign(x_i^k-x_j^k)- \E_{\bar\zeta}[\sum_{j \in \mathcal{R}_i(\bar\zeta)}\bar{s}_{ij}] \|^2}\nonumber\\
		\leq& \black{8\|\nabla\E_{\xi_i^k}[F(x_i^k,\xi_i^k)] +\nabla f_0(x_i^k) - \nabla \E_{\xi_i^k}[F(x_i^*,\xi_i^k)] -\nabla f_0(x_i^*)\|^2 + 32\lambda^2|\mathcal{R}_i|^2p.} \nonumber
	\end{align}
	Here the first inequality is due to $(a+b)^2\leq2a^2+2b^2$. To obtain the last inequality, we observe that \black{$sign(x_i^k-x_j^k),\bar{s}_{ij}\in[-1,1]^p$}. Combining \eqref{e3} and \eqref{e4}, we have
	\begin{align}\label{e4-tmp-tmp}
		&\black{\E_k}\|\nabla F(x_i^k,\xi_i^k)+ \lambda\sum_{j \in \mathcal{R}_i^k} sign(x_i^k-x_j^k)+\lambda\sum_{j \in \mathcal{B}_i^k}sign(x_i^k-z_j^k)+\nabla f_0(x_i^k)\|^2 \\
		\leq& \black{8\|\nabla\E_{\xi_i^k}[F(x_i^k,\xi_i^k)] +\nabla f_0(x_i^k) - \nabla \E_{\xi_i^k}[F(x_i^*,\xi_i^k)] -\nabla f_0(x_i^*)\|^2
			+ 32\lambda^2|\mathcal{R}_i|^2p + 4 \lambda^2 |\mathcal{B}_i|^2p + 2\delta_i^2.} \nonumber
	\end{align}

	For the third term at the right-hand side of \eqref{e1}, similar to the proof from \eqref{b2} to \eqref{b8}, we reach
	\begin{align}\label{b8-again}
		&-2\black{\E_k}\langle \nabla F(x_i^k,\xi_i^k)+\lambda\sum_{j \in \mathcal{R}_i^k} sign(x_i^k-x_j^k)+\nabla f_0(x_i^k), x_i^k-x_i^*\rangle\\
		\leq & -2\bigg\langle\lambda\sum_{j \in \mathcal{R}_i^k}sign(x_i^k-x_j^k)-\lambda\black{\E_{\bar\zeta}} \big[\sum_{j \in \mathcal{R}_i(\black{\bar\zeta})}\black{\bar{s}_{ij}}\big] , x_i^k-x_i^*\bigg\rangle \nonumber
		-\frac{2u_iL_i}{u_i+L_i} \|x_i^k-x_i^*\|^2 \nonumber \\
		&- \frac{2}{u_i+L_i} \ \black{\|\nabla \E_{\xi_i^k}[F(x_i^k,\xi_i^k)] + \nabla f_0(x_i^k) - \nabla \E_{\xi_i^k}[F(x_i^*,\xi_i^k)] - \nabla f_0(x_i^*)\|^2}. \nonumber
	\end{align}
	Notice it holds for any $\epsilon>0$ that
	\black{
		\begin{align}\label{b8-plus}
			&-2\bigg\langle\lambda\sum_{j \in \mathcal{R}_i^k}sign(x_i^k-x_j^k)-\lambda\E_{\bar\zeta}\big[\sum_{j \in \mathcal{R}_i(\bar\zeta)}\bar{s}_{ij}\big] , x_i^k-x_i^*\bigg\rangle \\
			\leq &\frac{\epsilon}{2} \|x_i^k-x_i^*\|^2 + \frac{2\lambda^2}{\epsilon} \| \sum_{j \in \mathcal{R}_i^k}sign(x_i^*-x_j^*)- \E_{\bar\zeta}\big[\sum_{j \in \mathcal{R}_i(\bar\zeta)}\bar{s}_{ij}\big] \|^2 \nonumber \\
			\leq &\frac{\epsilon}{2} \|x_i^k-x_i^*\|^2 + \frac{8\lambda^2 |\mathcal{R}_i|^2p}{\epsilon}, \nonumber
		\end{align}
		where the last inequality comes from that the candidate non-zero entries of $sign(x_i^k-x_j^k)$ and $\bar{s}_{ij}$ are bounded by $1$.} Note that due to the use of these loose upper bounds, the second term at the right-hand side of \eqref{b8-plus} is proportional to $|\mathcal{R}_i|^2$. As a consequence, the eventually derived size of convergence neighborhood is also monotonically increasing when $|\mathcal{R}_i|^2$ increases.
	%
	%
	Combining \eqref{b8-again} and \black{\eqref{b8-plus}} yields
	\begin{align}\label{b8-plusplusplus}
		&\black{\E_k}\langle \nabla F(x_i^k,\xi_i^k)+\lambda\sum_{j \in \mathcal{R}_i^k} sign(x_i^k-x_j^k)+\nabla f_0(x_i^k), x_i^k-x_i^*\rangle\\
		\leq&\frac{\epsilon}{2} \|x_i^k-x_i^*\|^2 + \frac{8\lambda^2 |\mathcal{R}_i|^2p}{\epsilon} - \frac{2u_iL_i}{u_i+L_i} \|x_i^k-x_i^*\|^2 \nonumber \\
		&-\frac{2}{u_i+L_i} \black{\|\nabla \E_{\xi_i^k}[F(x_i^k,\xi_i^k)] + \nabla f_0(x_i^k) - \nabla \E_{\xi_i^k}[F(x_i^*,\xi_i^k)] - \nabla f_0(x_i^*)\|^2}. \nonumber
	\end{align}

	For the last term at the right-hand side of \eqref{e1}, it holds for any $\epsilon>0$ that
	\begin{align}\label{b9-again}
		&-2\langle \lambda\sum_{j \in \mathcal{B}_i^k} sign(x_i^k-z_j^k), x_i^k-x_i^*   \rangle \\
		\leq& \frac{\epsilon}{2}\|x_i^k-x_i^*\|^2 + \frac{2\lambda^2}{\epsilon} \|\sum_{j \in \mathcal{B}_i^k} sign(x_i^k-z_j^k)\|^2 \nonumber\\
		\leq& \frac{\epsilon}{2} \|x_i^k-x_i^*\|^2 + \frac{2\lambda^2|\mathcal{B}_i|^2p}{\epsilon}. \nonumber
	\end{align}

	Substituting \eqref{e4-tmp-tmp}, \eqref{b8-plusplusplus} and \eqref{b9-again} into \eqref{e1}, we have:
	\black{
		\begin{align}\label{e10}
			&\E_k\|x_i^{k+1}-x_i^*\|^2\\
			\leq&\big(1 -\alpha^{k}(\frac{2u_iL_i}{u_i+L_i}-\epsilon)\big)\|x_i^k-x_i^*\|^2
			+(\alpha^{k})^2( 32\lambda^2|\mathcal{R}_i|^2p+4\lambda^2|\mathcal{B}_i|^2p+2\delta_i^2) \nonumber + \alpha^k (\frac{2 \lambda^2 |\mathcal{B}_i|^2 p}{\epsilon} + \frac{8 \lambda^2 |\mathcal{R}_i|^2 p}{\epsilon})\\
			&-2\alpha^{k} (\frac{1}{u_i+L_i}-4\alpha^{k}) \|\nabla \E_{\xi_i^k}[F(x_i^k,\xi_i^k)] + \nabla f_0(x_i^k) - \nabla \E_{\xi_i^k}[F(x_i^*,\xi_i^k)] - \nabla f_0(x_i^*)\|^2. \nonumber
		\end{align}
	}
	Due to the step size rule, we have \black{$\frac{1}{u_i+L_i}-4\alpha^{k} \geq 0$. Therefore, we drop the last term} at the right-hand side of \eqref{e10}. Also noticing the definition of $\eta$, we rewrite \eqref{e10} into
	\black{
		\begin{align}\label{e10-tmp}
			&\E_k\|x_i^{k+1}-x_i^*\|^2\\
			\leq&\big(1 -\eta\alpha^{k}\big)\|x_i^k-x_i^*\|^2
			+(\alpha^{k})^2( 32\lambda^2|\mathcal{R}_i|^2p+4\lambda^2|\mathcal{B}_i|^2p+2\delta_i^2) \nonumber + \alpha^k (\frac{2 \lambda^2 |\mathcal{B}_i|^2 p}{\epsilon} + \frac{8 \lambda^2 |\mathcal{R}_i|^2 p}{\epsilon}).
		\end{align}
	}
	\noindent\black{\textbf{Step 2.} With \eqref{e10-tmp}, the rest of the proof follows Step 3 of Theorem \ref{theorem2} by taking expectation over $\{x_i^l: l\le k, i \in \mathcal{R}\}$.}

	\section{Proof of Theorem \ref{theorem6}}
	\label{proof-theorem6}
	
	\noindent\textbf{Proof.} When $\lambda\geq\lambda_0$, combining Theorems \ref{theorem4} and \ref{theorem5}, we directly reach \eqref{th6-1}. When $0<\lambda<\lambda_0$, we have
	\begin{align}
		\E[\|x^{k+1}-[\tilde{x}^*]\|^2]\leq 2\E[\|x^{k+1}-x^*\|^2] +2 \E[\|x^*-[\tilde{x}^*]\|^2].
	\end{align}
	Because $\E\|x^{k+1}-x^*\|^2 \leq \frac{\Delta_6}{k+1}+\overline{\alpha}\Delta_5 \ $ and $\ \E[\|x^*-[\tilde{x}^*]\|^2]\leq\Delta_7$, we reach \eqref{th6-2} and complete the proof.

\end{document}